\documentclass[12pt,authoryear]{elsarticle}
\usepackage{endnotes}
\usepackage{amsfonts}
\usepackage{amsmath}
\usepackage{eurosym}
\usepackage{amssymb}
\usepackage{graphicx}%
\usepackage{caption}
\usepackage{subcaption}
\usepackage{multirow}
\usepackage{color}
\usepackage{float}
\usepackage{hyperref}
\hypersetup{colorlinks,allcolors=black}
\usepackage[margin=2.5cm]{geometry}
\usepackage{natbib}
\bibpunct[, ]{(}{)}{,}{a}{}{,}%
\newtheorem{theorem}{Theorem}

\usepackage{array}
\newcolumntype{+}{>{\global\let\currentrowstyle\relax}}
\newcolumntype{^}{>{\currentrowstyle}}
\newcommand{\rowstyle}[1]{\gdef\currentrowstyle{#1}%
	#1\ignorespaces
}

\journal{Energy Economics}
\begin{document}
	
\begin{frontmatter}
\title{Joint optimization of sales-mix and generation plan for a large electricity producer}
\author{Paolo Falbo\fnref{label2}}
\ead{paolo.falbo@unibs.it }
\address[label2]{Department of Economics and Management, University of Brescia, C.da S. Chiara n. 50, Brescia, Italy}

\author{Carlos Ruiz\corref{cor1}\fnref{label1,label3}}
\ead{carlos.ruiz@uc3m.es}
\address[label1]{Department of Statistics \& UC3M-BS Institute for Financial Big Data (IFiBiD), University Carlos III de Madrid, Avda. de la Universidad, 30,	28911-Legan\'es, Spain}


\cortext[cor1]{Corresponding author.}
\cortext[]{We would like to acknowledge the help of Maren Schmeck at the early stages of this work.}
\fntext[label3]{The author gratefully acknowledge the financial support from the Spanish government through projects PID2020-116694GB-I00 and from the Madrid Government (Comunidad de Madrid) under the Multiannual Agreement with UC3M in the line of ``Fostering Young Doctors Research'' (ZEROGASPAIN-CM-UC3M), and in the context of the V PRICIT (Regional Programme of Research and Technological Innovation.}

\begin{abstract}
The paper develops a typical management problem of a large power producer (i.e., he can partly influence the market price). In particular, he routinely needs to decide how much of his generation it is preferable to commit to fixed price bilateral contracts (e.g., futures) or to the spot market. However, he also needs to plan how to distribute the production across the different plants under his control. The two decisions, namely the sales-mix and the generation plan, naturally interact, since the opportunity to influence the spot price depends, among other things, by the amount of the energy that the producer directs on the spot market. We develop a risk management problem, since we consider an optimization problem combining a trade-off between expectation and conditional value at risk of the profit function of the producer. The sources of uncertainty are relatively large and encompass demand, renewables generation and the fuel costs of conventional plants. We also model endogenously the price of futures in a way reflecting an information advantage of a large power producer. In particular, it is assumed that the market forecast the price of futures in a naive way, namely not anticipating the impact of the large producer on the spot market. The paper provides a MILP formulation of the problem, and it analyzes the solution through a simulation based on Spanish power market data. 
\end{abstract}

\begin{keyword}
	CVaR \sep demand uncertainty \sep electricity industry\sep futures market \sep renewable uncertainy \sep risk aversion \sep spot market.
\end{keyword}

\end{frontmatter}
\definecolor{light-blue}{rgb}{0.50,0.55,1.00}%
\definecolor{dark-blue}{rgb}{0.15,0.15,.60}%

\section{Introduction}

Electricity generation can be sold in two ways mainly: on the spot market and
through fixed-price contracts (i.e., bilateral, forwards and futures). We will
refer in this paper to futures as for the typical fixed-price contract, but
the same arguments could apply to the other types. Deciding how to optimally
distribute power generation between the two channels is known as the
\textit{sales-mix} problem. Such decision affects the profit of the producer
significantly. Since futures contracts fix the sales price and the quantity of
energy to be delivered, they introduce a fixed component in the profit function.

However, fixed price contracts do not necessarily reduce the volatility of
profits. Spot prices tend to move according to direct generation costs, so
directing the output to the spot market turns out in some cases to reduce
rather than increase volatility of profits. A second major advantage of
selling generation on the spot market is that of using superior information.
Bids on the spot market are submitted knowing generation costs, the residual
demand and the capacities made available by competitors.

Such information is key to optimize the\textit{ generation plan}, that is how
different plants and their capacities are assigned to satisfy the demand. In
particular, directing the generation of cheapest technologies to satisfy
either futures contracts or the spot market can have a significant impact on
the spot price and consequently on profit. So, even if the sales-mix decision
is usually taken before that of the generation plan, the two are clearly
linked and so a major producer needs to optimize them jointly.

In this paper we model such a joint optimization problem, considering several
important general settings of the electricity markets. In the first place the
model accounts for the interplay between the spot and the futures prices of
electricity. In particular, the futures price is modeled as a risk neutral
expectation of the spot price that the market obtains by (partially)
anticipating the optimal behaviour of the major producer. Secondly, the model
includes the impact of the sales-mix decision on the spot demand. In
particular the demand of electricity arriving on the spot market is reduced by
the amount of energy that the major producer commits to the futures
contracts. Moreover, to account for the relevance of the renewable energy (RE)
technologies and risk aversion, the problem is formulated in four versions
combining the case where the major producer is endowed with high and low
percentage of renewables of his total capacity, and he is either risk neutral
or risk averse.

All these features make this work original with respect to the existing
literature. In particular, this works improves the findings of
\cite{FalboRuiz(2019)}, which also optimize the joint sales-mix and generation
plan decisions, but taking futures price as an exogenous parameter and letting
market demand independent of the fixed-price contracts signed by the major producer.

Modeling futures price as an endogenous variable allows this work to analyze
under which conditions positive or negative risk premia emerge in electricity
markets. Several other works have analyzed the impact of futures or forward
transactions. \cite{Allaz_1992} and \cite{Allaz_Villa_1993} study a model of
the relationship between futures and spot markets, in the contest of a duopoly
model, show that market prices may decrease with the introduction of futures
trading. That model has been also considered by \cite{Dong_Liu_2007} in the
context of a supply chain where players are risk averse. The impact of
introducing fixed-price contracts on supply chain efficiency has also been
studied by \cite{Mendelson_2007}. These authors show that supply and demand
uncertainty play a key role in the relationship between forward and spot
trading. According to \cite{Popescu_seshadri_2013} the uncertainty of demand
is also found to affect significantly the optimal sales policy in forward and
spot markets. In particular, the authors find that demand elasticity conditions
significantly forward and spot trading volume.

Another important feature of the analysis in this work is concerned with the
correlation between demand of electricity and renewable generation. This
element is often neglected in the literature of futures electricity markets,
yet it is a central factor affecting the expected prices of electricity.
\cite{FalboRuiz(2019)} report empirical evidence that such correlation can
change significantly across different regions and through time. They also show
how such factor influences the optimal policy of a power company. However,
their model did not allow them to check the impact on the futures prices. Some
papers have focused on the correlations among other key factors characterizing
the electricity markets. In their analysis on price forecasting and risk
management for an electricity producer \cite{Alvarez_2010} focus on time
varying correlations between electricity prices and demand. Same type of
correlation has also been analyzed in \cite{Boroumand_2015} to the purpose of
setting up an optimal hedging portfolio management for retailers seeking to
meet demand in different intraday time segments. In this work such correlation
is endogenous as it is the result of the interaction between demand and the
strategic bidding of the major producer.

The impact of renewables on the power sector has also been analyzed under
different aspects. \cite{Bell_2015} find a small correlation between
electricity demand and wind speed in Australia, tending to increase during
their observation period (2010-2012). They analyze the benefit of meeting
electricity demand in the different Australian states (without the use of
storages) resulting from this increasing correlation. Still focusing on the
Australian market \cite{Cutler2011} also find strong evidence of negative
correlation between wind generation and electricity price. In
\cite{Chaiamarit_2014} the authors assess the impact of several renewable
sources on the demand of electricity and, consequently, find meaningful
implication on operating costs. Large evidence of a negative impact of RES on
electricity prices in Germany has been found by \cite{Dillig2016}.

The correlation among the key variables of electricity markets is also central
to the problem of diversifying the technology portfolio (also referred as the
energy-mix problem). \cite{Adabi_2016} and \cite{Delarue2011} resort to the
theory of financial portfolios. In particular, the latter paper combines short
term generation and long term investment optimization. Interestingly such work
shows that wind capacity can help to lower the risk on generation cost.
However, that paper focuses just on the cost side the problem, so it misses the
correlation between electricity prices and producing factors (usually gas
cost), which drives the risk of electricity generators. Vice versa, a paper
explicitly considering the opportunity to reduce risk leveraging on the
correlation between electricity price and generation costs is that of
\cite{FalboEtAl2010}. The connection of this paper with this work is that in
both cases (part of) the solution of the optimization problem is set with the
sales-mix between spot market and bilateral contracts. However, in that work
the decision is taken by individual price-taker producer, while here a game
theoretic approach is considered between a producer and more competitors, so
that electricity price is an outcome depending on the solution adopted in the game.

Focusing on electricity markets, \cite{Niu_2005} propose a supply function
equilibrium model to derive the optimal bidding strategy of firms considering
forward contracts. An oligopolistic model is used by \cite{Anderson_Hu_2008}
to study how futures trading may increase social welfare when retailers have
market power. \cite{Bushnell_2007} extend the Allaz Vila model to multiple
firms and increasing costs to show that suppliers' market power plays a key
role in the interaction between forward and spot trading. A stochastic
programming framework is presented by \cite{Carrion_et_al_2007} to derive the
optimal futures and spot trading by a risk-averse and a price-taking retailer.
Similar to the proposed model, \cite{Conejo_et_al_2008} analyze the optimal
involvement of a price-taker electricity producer in the futures market
considering risk aversion and exogenous prices. An equilibrium model is
proposed by \cite{Aid_2011} to show that vertical integration plays a key role
in the relationship between spot and forward markets with demand uncertainty.
Another equilibrium model is presented in \cite{Ruiz_2012} to analyze the
equilibria in futures and spot markets with oligopolistic generators and
conjectural variations. Different types of futures contracts, and their impact
on the market outcomes, are tested in a game theoretical setting by
\cite{Oliveira_et_al_2013}. \cite{Fanzeres_2015} derive the optimal
contracting strategy of a price-taker and risk averse energy trading company
which operates with a renewable portfolio. Considering a medium-term model,
\cite{Mari_et_al_2017} propose a heuristic algorithm to derive the optimal
involvement in pool and bilateral contracts of a generation company.

This paper contributes to the state-of-the-art on this topic by presenting a
new complex and computationally efficient tool for decision making under
uncertainty. Power producers can benefit from it in order to decide which is
the optimal level of commitment in a futures market, together with the most
adequate strategy (generation plan) to supply that energy. Similarly, market
operators can use it to evaluate the impact of competition and risk aversion
in two-stage markets and infer how these may evolve in the coming years under
the smart grid paradigm. More specifically, the contributions of this work are seven-fold:

\begin{itemize}
	\item[1)] To derive the simultaneous optimal sales-mix and generation plan of
	a generation company with different types of technologies: nuclear, renewables
	and conventional generation units.
	
	\item[2)] To model the
	multivariate risk generated by the interaction of by demand
	load, renewable generation, and generation costs uncertainties.
	
	\item[3)] To account for the generators market power by the endogenous
	formation of spot prices.
	
	\item[4)] To formulate a bi-level stochastic optimization model to derive the
	optimal strategy in 1). The model is converted into an equivalent single level
	Mathematical Problem with Equilibrium Constraints (MPEC) that is latter
	linearized without approximation into an equivalent Mixed-Integer Linear
	Problem (MILP) easy to solve by available branch-and-cut solvers.
	
	\item[5)] To extend the model formulation in 4) to account for risk via the
	Conditional Value-at-Risk (CVaR).
	
	\item[6)] To analyze the impact of demand, renewable generation and direct costs volatility adopting a set of realistic numerical simulations based on Spanish Electricity Market data.
	
	\item[7)] To analyse the impact of the ``information gap hypothesis'' on the way spot price expectations are formed and on the prevalence of negative risk premia\footnote{In the financial literature of futures contracts, the risk premium is referred as the difference between the time $t$ price of a futures closing at $T$, $p^f_{t,T}$ and the expectation of the spot price of the underlying asset given the information available at time $t$, $\mathbb{E}[p_T^{s}|I_t]$.} in the futures markets of the power sector.
\end{itemize}

\section{Problem settings}

We consider the problem faced by a strategic electricity producers that participates in a two-stage electricity market composed of a futures market that is followed by a spot market. The decisions taken in the futures market have a direct impact on the spot market, as they condition both the residual demand in that market together with the available production capacities for each technology.  Index $i=1,\dots,I$ is used to refer to the different generating units owned by the strategic producer while index $j=1,\dots,J$ characterize each generating unit own by rival producers\footnote{From the strategic producer's point of view, it is indifferent to assume a single or several rival producers owning the generating units $j=1,\dots,J$.}. These units $i$ or $j$ may be based on different generating technologies (renewable or based on fossil fuels).

We generalize the model to consider a risk averse player that needs to take two type of strategic decisions: i) how much to trade through futures and hence condition his remaining generation capacity and residual demand in the spot market, and ii) how to distribute the generation between its different units $i$ to satisfy both the futures and spot commitments.

Regarding the futures market price $p^f$ formation, we generalize the classical non-arbitrage condition, which states that this price equals the expected spot price. In this model, we propose to assume that the resulting futures price $p^f$ is formed based on ``subjective''  expectations that the different players would have on the subsequent spot market, which may be accurate or not. This will allow us testing different levels of players' ``expertise'' in anticipating the spot market, and to evaluate whether the market can show positive or negative risk premia.

\subsection{Bilevel Stochastic Formulation}
The decision-making process of the strategic produce is mathematically modeled as a bilevel stochastic problem. The upper-level problem characterizes the maximization of the strategic producer profit while anticipating two lower-level problems that reproduce the clearing of the spot market: one of them is used to compute the futures market price $p^f$ (based on imperfect information about the spot market) and the other one is used to anticipate the true dispatch of its generating units (and resulting spot price). Uncertainty is incorporated to the model via scenarios (index $\omega=1,\dots,\Omega$), reproducing different potential realizations of the uncertain parameters (renewable capacity, fuel costs, total demand, etc.) and hence different clearing of the spot market.

The strategic produces faces the following bi-level optimization problem:
\begin{subequations}\label{bilevel}
	\begin{align}
		&  \max_{q^{ft},p^f,q_{i\omega}^{f},q_{i\omega}^{s},p^{s}_{\omega},\breve{p}^{s}_{\omega},\mathcal{Y}}\ \mathbb{H} \left[  p^{f}q^{ft}+p_{\omega}^{s}\sum^{I}_{i=1}q_{i\omega
		}^{s}-\sum^{I}_{i=1}c^{s}_{i\omega}(q_{i\omega}^{f}+q_{i\omega} ^{s}) \quad \forall \omega\right]\label{bi_of}
		\\
		&  \mbox{s.t.}\nonumber\\
		&  \quad q^{ft}=\sum^{I}_{i=1} q_{i\omega}^{f}\qquad\forall\omega\label{bi_cons_1}\\
		&  \quad 0\leq q_{i\omega}^{f}\leq Q_{i\omega}^{s}\quad\forall i,\forall
		\omega \label{bi_cons_2}\\
		& \quad p^f = \sum_{\omega=1}^{\Omega} \sigma_{\omega} \breve{p}^{s}_{\omega}\label{bi_cons_3}\\
		&  \quad \breve{p}^s_{\omega}\in\left\{ \mbox{\textit{``Naive'' Spot Market Clearing} (\ref{naive})}\right\}\quad \forall \omega \label{bib_cons_4}\\
		&  \quad q_{i\omega}^{s},p^{s}_{\omega}\in\left\{\mbox{\textit{``Strategic'' Spot Market Clearing} (\ref{strategic})}\right\}\quad \forall \omega \label{bib_cons_5}
	\end{align}
\end{subequations}
where $\mathbb{H}[\hphantom{x}]$ represents a general risk measure computed over the uncertain profit (expected value or value at risk for instance). Hence, the expression between brackets in (\ref{bi_of}) represents the producers profit per scenario $\omega$. The first term stands for the revenue obtained in the futures market, i.e., futures price $p^f$ times futures total production $q^{ft}$. Note that this revenue is scenario independent, as it depends on first-stage (here and now) decisions. The second term in (\ref{bi_of}) represents the total income from selling each generating unit's spot production quantity $q_{i\omega}^s$ at the spot price $p^s_{\omega}$, for each scenario $\omega$. Finally, the last term accounts for the total generation costs, where $c_{i\omega}^s$ is the variable (linear) generation costs for each unit $i$ and scenario $\omega$. Observe that the futures market entail a physical delivery of electricity and hence contributes to the generation costs, where $q_{i\omega}^f$ is the amount of energy generated by unit $i$ to satisfy the futures contract. It should be emphasized that despite the total futures trading $q^{ft}$ is scenario independent (first-stage decision), this is not the case of how this production is latter distributed among the generating units ($q_{i\omega}^f$ depends on $\omega$), decision that can be postponed to the second stage (the futures energy delivery takes place at the same time that the spot market is cleared). Indeed, constraint (\ref{bi_cons_1}) ensures that the total quantity traded in the futures $q^{ft}$ is satisfied for all the scenarios $\omega$, by different mix of generating units.

Constraints (\ref{bi_cons_2}) enforce, for each scenario $\omega$ and generating unit $i$, that the energy quantities generated to supply the futures contract is below the maximum generation capacity for each unit $i$, i.e., $Q_{i\omega}^s$. These capacities may be uncertain (depend on scenario $\omega$), specially if associated to a renewable energy source.

Furthermore, we assume that the futures price $p^f$ represents the expectation of the spot price under partial information. Since we discretize the probability space into scenarios, this is computed as shown in (\ref{bi_cons_3}), where $\sigma_{\omega}$ is the probability associated to each scenario and $\breve{p}^{s}_{\omega}$ is the corresponding estimated spot price. In particular, $\breve{p}^{s}_{\omega}$ is obtained from the solution of the lower-level problem (\ref{bib_cons_4}) that, under some partial information, reproduces the spot market clearing for each scenario $\omega$. On the contrary, the ``true'' spot price and dispatched quantities per scenario, anticipated by the strategic producer ($p^{s}_{\omega}$ and $q_{i\omega}^{s}$, respectively), are computed in a lower-level problem that reproduces the actual clearing of the spot market (\ref{bib_cons_5}). The next subsections present the specific formulations of the lower-level problems (\ref{bib_cons_4}) and (\ref{bib_cons_5}), and how they are directly affected by the upper-decision variables.

\subsubsection{``Naive'' Spot Market Clearing}
Problem (\ref{naive}) below, reproduces the formation of the spot prices for each scenario $\omega$, assuming partial information about the decisions taken by the strategic producer. To differentiate them with the actual spot market outcomes, both the price and the dispatch quantities are marked with the accent ``$\breve{\hphantom{x}}$''.
\begin{subequations}\label{naive}
	\begin{align}
		&  \min_{\breve{q}_{i\omega}^{s},\breve{q}_{j\omega}^{sr}}\ \sum^{I}_{i=1} c_{i\omega}^{s}
		\breve{q}_{i\omega}^{s}+\sum^{J}_{j=1}c_{j\omega}^{sr}\breve{q}_{j\omega}^{sr} \label{naive_of}\\
		&  \mbox{s.t.}\nonumber\\
		&  \quad \sum^{I}_{i=1}
		\breve{q}_{i\omega}^{s}+\sum^{J}_{j=1}\breve{q}_{j\omega}^{sr}=D_{\omega}-\alpha q^{ft}\quad(\breve{p}^{s}_{\omega})\label{naive_bal}\\
		& \quad 0\leq \breve{q}_{i\omega}^{s}\leq Q_{i\omega}^{s}-\frac{Q^s_{i\omega}}{\sum_{i=1}^I Q^s_{i\omega}}\alpha q^{ft}\qquad\forall i \label{naive_cap}\\
		&  \quad  0\leq \breve{q}_{j}^{sr}\leq Q_{j\omega}^{sr}\label{naive_end}
		\qquad\forall j
	\end{align}
\end{subequations}
The objective function (\ref{naive_of}) represents the minimization of the total generating costs, performed by the ISO, to satisfy the demand. In particular the first term is the costs incurred by the units of the strategic producer, and the second term corresponds to the rival generators, i.e., costs $c_{j\omega}^{sr}$ times generation quantities $\breve{q}_{j\omega}^{sr}$. The balance constraint is imposed by (\ref{naive_bal}), where the total generation by strategic ($\sum^{I}_{i=1}\breve{q}_{i\omega}^{s}$) and non-strategic ($\sum^{J}_{j=1}\breve{q}_{j\omega}^{sr}$) units equals the residual demand in the spot market. The residual demand is composed of the total load in the system at scenario $\omega$ ($D_{\omega}$) minus the energy that is going to be supplied by the futures contract ($q^{ft}$), where $\alpha$ is a coefficient representing the correctness of the guess of rivals about the futures commitment of the strategic producer. More precisely, $\alpha=1$ represents a perfect guess, while values of $\alpha > 1$ or $\alpha < 1$ represent respectively over or under estimations. In this work we will assume $\alpha=1$.
The dual variable associated to this constraint, i.e., $\breve{p}^{s}_{\omega}$, can be interpreted as the spot price as, for each scenario $\omega$, it identifies the marginal costs of supplying an additional unit of residual demand. Indeed, this is the price that is used in (\ref{bi_cons_3}) to estimate the futures price.

The minimum and maximum production capacity for each strategic unit is limited by (\ref{naive_cap}). In particular, we assume that the maximum capacity equals the total available capacity for that unit ($Q_{i\omega}^s$) minus the share that it is used to supply the futures contract. Moreover, to characterize this ``naive'' vision of the spot market, we consider that the market players (both strategic and non-strategic) assume that the futures contract is supplied by each technology proportionally to their available capacities, i.e., actual capacity of a unit $i$ ($Q^s_{i\omega}$) over the total ($\sum_{i=1}^I Q^s_{i\omega}$). This assumption is relaxed in the second lower level problem (\ref{strategic}) to allow for strategic behaviors. Finally, constraint (\ref{naive_end}) imposes lower and upper ($Q^{sr}_{j\omega}$) bounds for the generation capacity of the non-strategic units.

\subsubsection{``Strategic'' Spot Market Clearing}

The spot market clearing problem that it is actually anticipated by the strategic producer is presented in (\ref{strategic}).
\begin{subequations}\label{strategic}
	\begin{align}
		&    \min_{q_{i\omega}^{s},q_{j\omega}^{sr}} \  \sum^{I}_{i=1}
		c^{s}_{i\omega}q_{i\omega}^{s}+\sum^{J}_{j=1}c^{sr}_{j\omega} q_{j\omega}%
		^{sr} \label{str_of}\\
		&  \mbox{s.t.}\nonumber\\
		&  \quad \sum^{I}_{i=1}q_{i\omega}^{s} +\sum^{J}_{j=1}q_{j\omega}^{sr}=D_{\omega}- q^{ft}\quad(p^{s}_{\omega
		})\label{str_bal}\\
		&  \quad 0\leq q_{i\omega}^{s}\leq Q_{i\omega}%
		^{s}-q^{f}_{i\omega}\qquad(\mu_{i\omega}^{\min},\mu_{i\omega}^{\max}%
		)\qquad\forall i\label{str_cap}\\
		&  \quad 0\leq q_{j\omega}^{sr}\leq Q_{j\omega}%
		^{sr}\qquad(\nu_{j\omega}^{\min},\nu_{j\omega}^{\max})\qquad\forall j\qquad \quad \label{str_end}
	\end{align}
\end{subequations}

Similar to (\ref{naive}), the objective function (\ref{str_of}) aims to minimize the total generating costs considering both the strategic (first term) and non-strategic (second term) units. Equation ($\ref{str_bal}$) represents the spot market balance where total generation (left term) matches residual demand (right term). Again, the dual variable associated with this constraint ($p^s_{\omega}$) represents the spot price per scenario $\omega$. Constraint (\ref{str_cap}) imposes the minimum and maximum generation capacity of unit $i$. As opposed to (\ref{naive_cap}), we assume that the strategic producer can decide, at the time that the spot market takes place (second stage), which amount of each unit's capacity is assigned to supply the futures and the spot market commitments. In this sense, $q_{i\omega}^f$ for every unit $i$ can be decided strategically to alter the marginal cost, and hence the price, of the spot market. Constraint (\ref{str_end}) sets lower and upper limits for the generation capacity of the non-strategic units. Dual variables are indicated at their corresponding constraints between parentheses.

\subsection{Solution Methodology}\label{sec:solution_methodology}
Bilevel problem (\ref{bilevel}) includes two lower level problems that complicate to find its global solution, either analytically or numerically. Apart from heuristic frameworks that cannot guarantee global optimality, most of standard approaches are based on replacing the lower-level problems by their first order optimality KKT (Karush-Kuhn-Tucher) conditions, and try to linearize the resulting single level problem. This approach is not completely useful in our case as there are still some nonlinear terms that cannot be linearized.

However, due to the particular two-stage decision structure of our problem (\ref{bilevel}), we propose an alternative exact solution methodology. The idea is to parameterize problem (\ref{bilevel}) in the first-stage decision variable $q^{ft}$, and explore its range of possible values. This is a convenient choice as $q^{ft}$ is a single variable whose entire feasible region can be easily explore. In particular, note that from (\ref{bi_cons_1}) and (\ref{bi_cons_2}) we can derive that $0\leq q^{ft}\leq \max_{\omega}\sum_{i=1}^{I}Q^s_{i\omega}$. The solution methodology works as follow:

\begin{itemize}
	\item[i)] Fix $q^{ft}$ to a previously-not-chosen value within the range $0\leq q^{ft}\leq \max_{\omega}\sum_{i=1}^{I}Q^s_{i\omega}$.
	\item[ii)]  If $q^{ft}$ is fixed then (\ref{naive}) can be solved as a linear problem and $\breve{p}^s_\omega$ can easily computed from the dual solution. Hence, the futures price $p^f$ can also be computed from (\ref{bi_cons_3}).
	\item[iii)] Then, by knowing $q^{ft}$ and $p^f$, the bilevel problem formed by (\ref{bi_of})-(\ref{bi_cons_3}) and (\ref{bib_cons_5}) can be recast as an MPEC by replacing problem (\ref{strategic}) by its KKT conditions. This results in (\ref{MPEC}) which can be linearized and solve efficiently as indicated in Section \ref{sec:linearization}.
	\item[iv)] Store the solutions from iii) and go back to i). Repeat until the entire range in i) is explored.
\end{itemize}

This procedure allows exploring the impact, in terms of profit distributions (one potential profit for each scenario $\omega$), of each possible value of the first stage decision $q^{ft}$. Then, by defining different risk-averse profit criteria, the optimal value of $q^{ft}$ can be selected, from risk-neutral to risk-averse.

Moreover, the proposed procedure is computationally efficient as, once the first-stage decision $q^{ft}$ is fixed and $p^f$ is obtained in ii), the computation of the solution of the MPEC in iii) can be decomposed per scenario $\omega$ (there are no linking variables) and hence parallelized in several machines.

\subsection{MPEC problem}\label{sec:MPEC}
The formulation of the MPEC problem derived from (\ref{bi_of})-(\ref{bi_cons_3}) and (\ref{bib_cons_5}) is presented in (\ref{MPEC}) below:
\begin{subequations}
	\label{MPEC}%
	\begin{align}
		&  \max_{q_{i\omega}^{f},q_{i\omega}^{s},p^{s}_{\omega},\mathcal{Y}}\ p^{f}q^{ft}+ \mathbb{H} \left[  p_{\omega}^{s}\sum^{I}_{i=1}q_{i\omega
		}^{s}-\sum^{I}_{i=1}c^{s}_{i\omega}(q_{i\omega}^{f}+q_{i\omega} ^{s}) \quad \forall \omega\right]\label{MPEC_of}\\
		&  \mbox{s.t.}\nonumber\\
		&  \quad q^{ft}=\sum^{I}_{i=1} q_{i\omega}^{f}\qquad\forall\omega
		\label{MPEC1}\\
		&  \quad0\leq q_{i\omega}^{f}\leq Q_{i\omega}^{s}\quad\forall i,\forall
		\omega\label{MPEC2}\\
		&\quad \sum^{I}_{i=1}q_{i\omega}^{s} +\sum^{J}_{j=1}q_{j\omega}^{sr}=D_{\omega}- q^{ft}\quad \forall \omega \label{MPEC3}\\
		&  \quad c^{s}_{i\omega}-p_{\omega}^{s}+\mu_{i\omega}^{\max}-\mu_{i\omega
		}^{\min} = 0\quad\forall i,\forall\omega\label{KKT1}\\
		&  \quad c^{sr}_{j\omega}-p_{\omega}^{s}+\nu_{j\omega}^{\max}-\nu_{j\omega
		}^{\min} = 0\quad\forall j,\forall\omega\label{KKT2}\\
		&  \quad0 \leq Q_{i\omega}^{s}-q_{i\omega}^{s}-q_{i\omega}^{f}\perp
		\mu_{i\omega}^{\max} \geq0 \qquad\forall i,\forall\omega\label{compl_1}\\
		&  \quad0 \leq q_{i\omega}^{s}\perp\mu_{i\omega}^{\min} \geq0 \qquad\forall
		i,\forall\omega\label{compl_2}\\
		&  \quad0 \leq Q_{j\omega}^{sr}-q_{j\omega}^{sr}\perp\nu_{j\omega}^{\max}
		\geq0 \qquad\forall j,\forall\omega\label{compl_3}\\
		&  \quad0 \leq q_{j\omega}^{sr}\perp\nu_{j\omega}^{\min} \geq0 \qquad\forall
		i,\forall\omega\label{compl_4}
	\end{align}
\end{subequations}
As $p^{f}$ and $q^{ft}$ are considered fixed at this stage, their product can be computed outside the loss function  $\mathbb{H}[\hphantom{x}]$ in (\ref{MPEC_of}). Constraints (\ref{MPEC3})-(\ref{compl_4}) represent the KKT conditions of problem (\ref{strategic}). In particular, (\ref{KKT1}) and (\ref{KKT2}) are the stationarity conditions, while (\ref{compl_1})-(\ref{compl_2}) are the complementary constraints (symbol $\perp$ indicates complementarity, i.e., $a\perp b \Leftrightarrow ab=0$).

\subsection{Linearization of problem MPEC problem (\ref{MPEC})}\label{sec:linearization}
Problem (\ref{MPEC}) can be exactly linearized so that it can rewritten as a MILP (Mixed-Integer-Linear Problem). These problems can be solved efficiently, with global optimality guarantees, with standard solvers. There are two sources of nonlinearities in problem (\ref{MPEC}) which can be tackled as shown below:

\begin{itemize}
	\item The complementarity constraints (\ref{compl_1})-(\ref{compl_4}) can be linearized with the Fortuny-Amat and McCarl linearization approach as follows.
	\begin{subequations}
		\label{lin_compl}%
		\begin{align}
			&  0 \leq Q_{i\omega}^{s}-q_{i\omega}^{s}-q_{i\omega}^{f}\leq u_{i\omega
			}^{\max}M\qquad\forall i,\forall\omega\\
			&  0 \leq\mu_{i\omega}^{\max}\leq(1-u_{i\omega}^{\max})M\qquad\forall
			i,\forall\omega\\
			&  0\leq q_{i\omega}^{s}\leq u^{\min}_{i\omega}M \qquad\forall i,\forall
			\omega\\
			&  0\leq\mu_{i\omega}^{\min} \leq(1-u^{\min}_{i\omega})M \qquad\forall
			i,\forall\omega\\
			&  0 \leq Q_{j\omega}^{sr}-q_{j\omega}^{sr}\leq v_{j\omega}^{\max}%
			M\qquad\forall j, \forall\omega\\
			&  0 \leq\nu_{j\omega}^{\max}\leq(1-v_{j\omega}^{\max})M\qquad\forall j,
			\forall\omega\\
			&  0\leq q_{j\omega}^{sr}\leq v^{\min}_{j\omega} M \qquad\forall
			j,\forall\omega\\
			&  0\leq\nu_{j\omega}^{\min} \leq(1-v^{\min}_{j\omega})M \qquad\forall
			j,\forall\omega\\
			&  \left\{  u_{i\omega}^{\max}, u_{i\omega}^{\min},v_{j\omega}^{\max
			},v_{j\omega}^{\min}\right\}  \in\{0,1\}
			\qquad\forall i, \forall j, \forall\omega.
		\end{align}
	\end{subequations}
	where $M$ is a sufficiently large constant.
	
	\item The nonlinear product $p_{\omega}^{s}\sum^{I}_{i=1}q_{i\omega
	}^{s}$ in (\ref{MPEC_of}). Multiplying the terms in (\ref{MPEC3}) by $p^s_{\omega}$ renders
	\begin{equation}
		p_{\omega}^{s}\sum^{I}_{i=1}q_{i\omega}^{s} = - p_{\omega}^{s}\sum^{J}_{j=1}q_{j\omega
		}^{sr} + p^s_\omega (D_{\omega}- q^f) \qquad\forall\omega. \label{lin_1}
	\end{equation}
	
	Then, by replacing $p^s_\omega$ from (\ref{KKT2}) in the first term of the right hand-side of (\ref{lin_1}) we obtain
	\begin{equation}
		p_{\omega}^{s}\sum^{I}_{i=1}q_{i\omega}^{s} = -\sum_{j=1}^J(c_{j\omega}+\nu_{j\omega}^{\max}-\nu_{j\omega}^{\min})q_{j\omega}^{sr}+p^s_\omega (D_{\omega}- q^f) \qquad\forall\omega,\label{lin_2}
	\end{equation}
	which combined with the complementarity constraints (\ref{compl_3}) and (\ref{compl_4}) results in the following linear expression
	\begin{equation}
		p_{\omega}^{s}\sum^{I}_{i=1}q_{i\omega}^{s} = - \sum^{J}_{j=1}c^{sr}_{j\omega}q_{j\omega}^{sr} -\sum
		^{J}_{j=1}\nu_{j\omega}^{\max}Q_{j\omega}^{sr}+p^s_\omega (D_{\omega}- q^f) \qquad\forall\omega, \label{linear3}
	\end{equation}
	which is a linear formulation of the bilinear products $p_{\omega}^{s}\sum
	^{I}_{i=1}q_{i\omega}^{s}$.
\end{itemize}

\section{Numerical application}
\subsection{Data}
The considered power system configuration is based on data from \cite{FalboRuiz(2019)}, which reproduces the electricity mix from the Spanish Power System \citep{ESIOS_2018}. Hence, we model a power system that includes five possible electricity generation technologies: nuclear-, wind-, solar-, coal- and gas-based. To account for a strategic producer with sufficient market power, we assume that 25\% of the total system power capacity is owned by him. This capacity is provided by an energy mix that includes $I=8$ generating units of different types (technologies). The remainder 75 \% of the system capacity belongs to the rival producers and is gathered in $J=8$ units of different generation technologies.

The summary of the considered data set is presented in Table \ref{table_data}. It includes the unit index ($i$ for strategic and $j$ for rival producer) (first column), the marginal costs (second column), the generation capacity (third column) and the cumulative capacity (fourth column) for each generating unit. Units are ordered from lowest to highest marginal costs, so that the second and fourth columns can be considered to build the aggregated supply curve of the system. The costs associated to conventional units (coal and gas) and the available capacities of renewable (solar and wind), are considered uncertain. Hence, these are summarized by their mean ($\mu$) and standard deviation ($\sigma$) values. Similarly, we consider that the system faces an inelastic and uncertain demand characterized by an expected value $\mu=19000 MWh$ and a standard deviation $\sigma =3800$. We can anticipate that, in the absence of futures trading, the intersection of the aggregate supply curve and the inelastic demand may result in a market price determined by the marginal costs of gas or coal units. Nevertheless, the marginal technology may change once uncertainty is introduced by in terms of possible scenario realizations.

\begin{table}[ht]\centering
	\caption{Generating units cost and capacity data}\smallskip
	\resizebox{8cm}{!}{
		\begin{tabular}{l c c c c }
			Gener. & Type & Marg. Costs  & Capacity & Exp. Cumul.     \\ 
			Units  & Tech. & (\euro/$MWh$) &  ($MWh$) & Cap. ($GWh$)    \\ \hline
			\multirow{ 2}{*}{$i=1$} & \multirow{ 2}{*}{nuclear}& $\mu=1\text{e-}3$ &\multirow{ 2}{*}{1250} & \multirow{ 2}{*}{1.25} \\ 
			&  & $\sigma=1\text{e-}6$ & &  \\ \hline
			\multirow{ 2}{*}{$j=1$} & \multirow{ 2}{*}{nuclear}& $\mu=1\text{e-}3$ &\multirow{ 2}{*}{3750} & \multirow{ 2}{*}{5.00} \\ 
			&  & $\sigma=1\text{e-}6$ & &  \\ \hline
			\multirow{ 2}{*}{$i=2$} & \multirow{ 2}{*}{wind}& $\mu=1\text{e-}3$ &$\mu=692.64$ & \multirow{ 2}{*}{5.69} \\ 
			&  & $\sigma=1\text{e-}6$ & $\sigma=207.792$ &  \\ \hline
			\multirow{2}{*}{$i=3$}&\multirow{2}{*}{wind}& $\mu=1\text{e-}3$ &$\mu=692.64$ & \multirow{ 2}{*}{6.38} \\ 
			&  & $\sigma=1\text{e-}6$& $\sigma=207.792$ & \\ \hline
			\multirow{ 2}{*}{$j=2$}&\multirow{ 2}{*}{wind}&$\mu=1\text{e-}3$ &$\mu=2077.92$ & \multirow{ 2}{*}{8.46} \\ 
			&  & $\sigma=1\text{e-}6$& $\sigma=623.376$  &     \\ \hline
			\multirow{ 2}{*}{$j=3$}&\multirow{ 2}{*}{wind}&$\mu=1\text{e-}3$ &$\mu=2077.92$ & \multirow{ 2}{*}{10.54} \\ 
			&  & $\sigma=1\text{e-}6$& $\sigma=623.376$  &  \\ \hline
			\multirow{ 2}{*}{$i=4$}&\multirow{ 2}{*}{solar}&$\mu=1\text{e-}3$ &$\mu=353.07$ & \multirow{ 2}{*}{10.89} \\ 
			&  & $\sigma=1\text{e-}6$& $\sigma=105.921$  & \\ \hline
			\multirow{ 2}{*}{$j=4$ }&\multirow{ 2}{*}{solar}& $\mu=1\text{e-}3$ &$\mu=1059.22$ & \multirow{ 2}{*}{11.95} \\ 
			&  & $\sigma=1\text{e-}6$& $\sigma=317.766$  &  \\ \hline
			\multirow{ 2}{*}{$j=5$}&\multirow{ 2}{*}{coal}&$\mu=32.99$ & \multirow{ 2}{*}{1500} & \multirow{ 2}{*}{13.45} \\ 
			&  & $\sigma=1.649$ &    &    \\ \hline
			\multirow{ 2}{*}{$i=5$}&\multirow{ 2}{*}{coal}&$\mu=36.64$ & \multirow{ 2}{*}{500} & \multirow{ 2}{*}{13.95}  \\ 
			&  & $\sigma=1.832$ &    &     \\ \hline
			\multirow{ 2}{*}{$j=6$}&\multirow{ 2}{*}{coal}&$\mu=39.78$ & \multirow{ 2}{*}{1500} & \multirow{ 2}{*}{15.45}   \\ 
			&  & $\sigma=1.989$ &    &      \\ \hline
			\multirow{ 2}{*}{$i=6$}&\multirow{ 2}{*}{coal}&$\mu=41.67$ & \multirow{ 2}{*}{500} & \multirow{ 2}{*}{15.95}  \\ 
			&  & $\sigma=2.083$ &    &     \\ \hline
			\multirow{ 2}{*}{$j=7$}&\multirow{ 2}{*}{gas}&$\mu=43.43$ & \multirow{ 2}{*}{9000} & \multirow{ 2}{*}{24.95} \\ 
			&  & $\sigma=2.171$ &    &     \\ \hline
			\multirow{ 2}{*}{$i=7$}&\multirow{ 2}{*}{gas}&$\mu=45.44$ & \multirow{ 2}{*}{3000} & \multirow{ 2}{*}{27.95} \\ 
			&  & $\sigma=2.272$ &    &    \\ \hline
			\multirow{ 2}{*}{$j=8$}&\multirow{ 2}{*}{gas}&$\mu=48.83$ & \multirow{ 2}{*}{9000} & \multirow{ 2}{*}{36.95}   \\ 
			&  & $\sigma=2.441$ &    &     \\ \hline
			\multirow{ 2}{*}{$i=8$}&\multirow{ 2}{*}{gas}&$\mu=56.68$ & \multirow{ 2}{*}{3000} & \multirow{ 2}{*}{39.95}   \\ 
			&  & $\sigma=2.834$ &    &     \\ \hline
	\end{tabular}} \label{table_data}
\end{table}

\subsection{Model application}

The application of the model has required first to generate a large set of
scenarios $\omega =1,2,...,300$ , each consisting of a realization of a
vector of random variables, namely:
				\begin{itemize}
					\item[\textbullet] Total electricity demand ($MWh$): ${D}_{\omega}$ ($D\sim N(19000,3800)$).
					
					\item[\textbullet] Strategic/rivals RES (solar and wind) generation capacity ($MWh$): $Q_{i\omega}^{s}$/$Q_{j\omega}^{s}$ for $i/j=2,3,4$.
					
					\item[\textbullet] Strategic/rivals marginal generating cost (\euro /$MWh$): $c_{i\omega}^{s}$/$c_{j\omega}^{s}$%
					for coal-based ($i/j=5,6$) and gas-based technologies ($i/j=7,8$).
				\end{itemize}
For simplicity, the random scenarios are generated based on Normal probability distributions, and no correlations are considered between the above random variables except those between wind power unit, i.e., the correlation between $Q^s_{i/j\omega}$ for $i/j=2,3$ is one. An in-depth analysis on potential impacts of correlations can be found in \cite{FalboEtAl2010} and \cite{FalboRuiz(2019)}. 

Next we fix a grid of the possible amounts $q^{f}=0,125,250,...,3000$ (MWh)
of energy committed with futures. For each value of $q^{ft}$, the optimization model is solved under $\Omega=300$ scenarios, as described in Section \ref{sec:solution_methodology}. The collection of the optimal solutions ($q_{i\omega}^{s},q_{i\omega}^{ft}$,...) will be stored to latter analyzed several properties of the strategic behavior of the producer (scenario dependent). Moreover, this methodology allows studying how the profit distribution evolves within the grid of values of $q^{ft}$, which is an insight not available if $q^{ft}$ was considered endogenous. These results are presented in the following section.

\subsection{Results}

As indicated, our model allows to analyze several interesting aspects of the behavior of
the strategic producer as well as to obtain several normative results for the power market. We comment them through several points.

\subsubsection*{\textbf{Futures reduce the risk of profits}}
\begin{figure}[H]
	\begin{subfigure}[b]{0.5\textwidth}
    	\includegraphics[height=5.5cm]{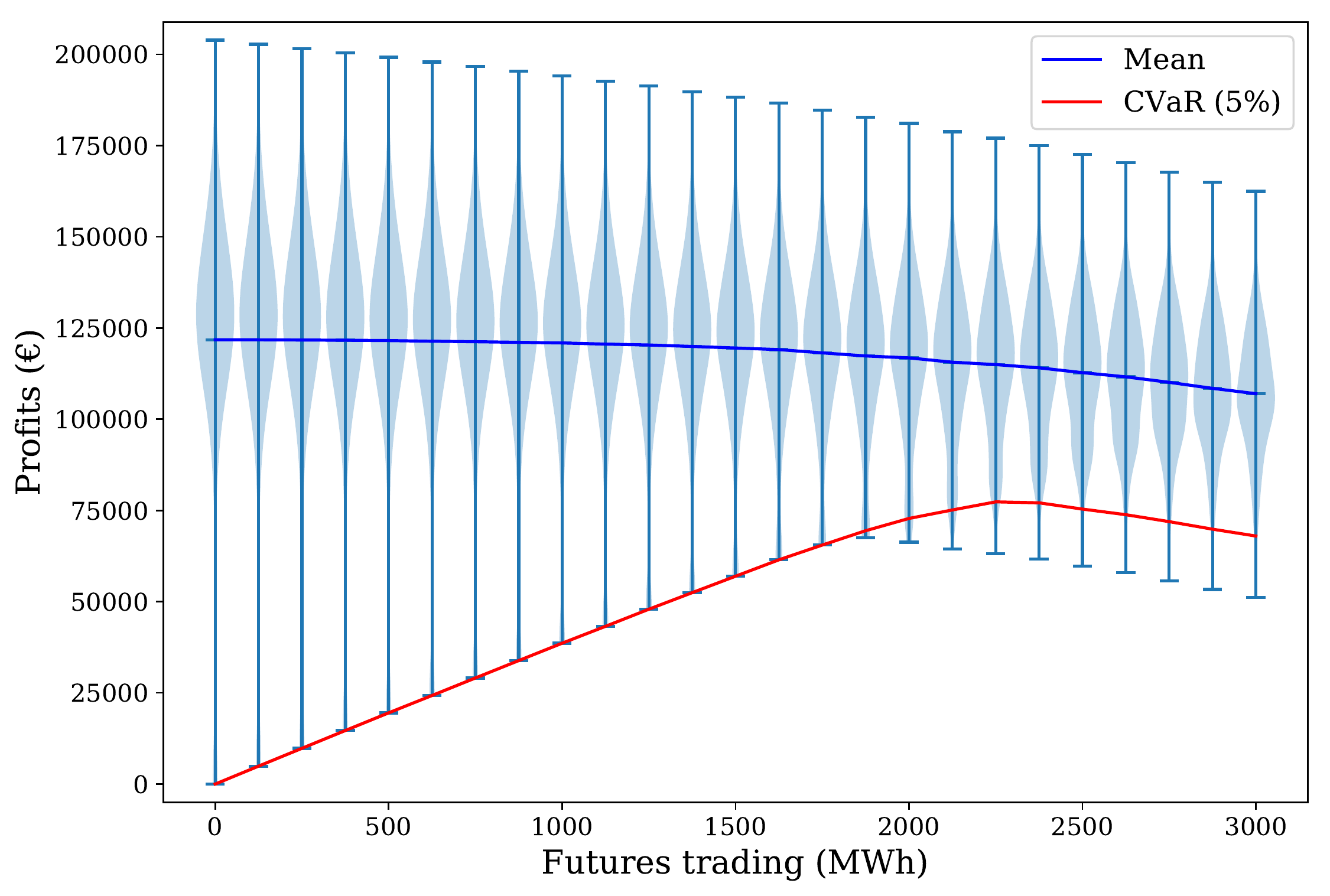}
    	\caption{Total profit}
    	\label{fig:profit_dist}
	\end{subfigure}	
	\begin{subfigure}[b]{0.5\textwidth}
    	\includegraphics[height=5.5cm]{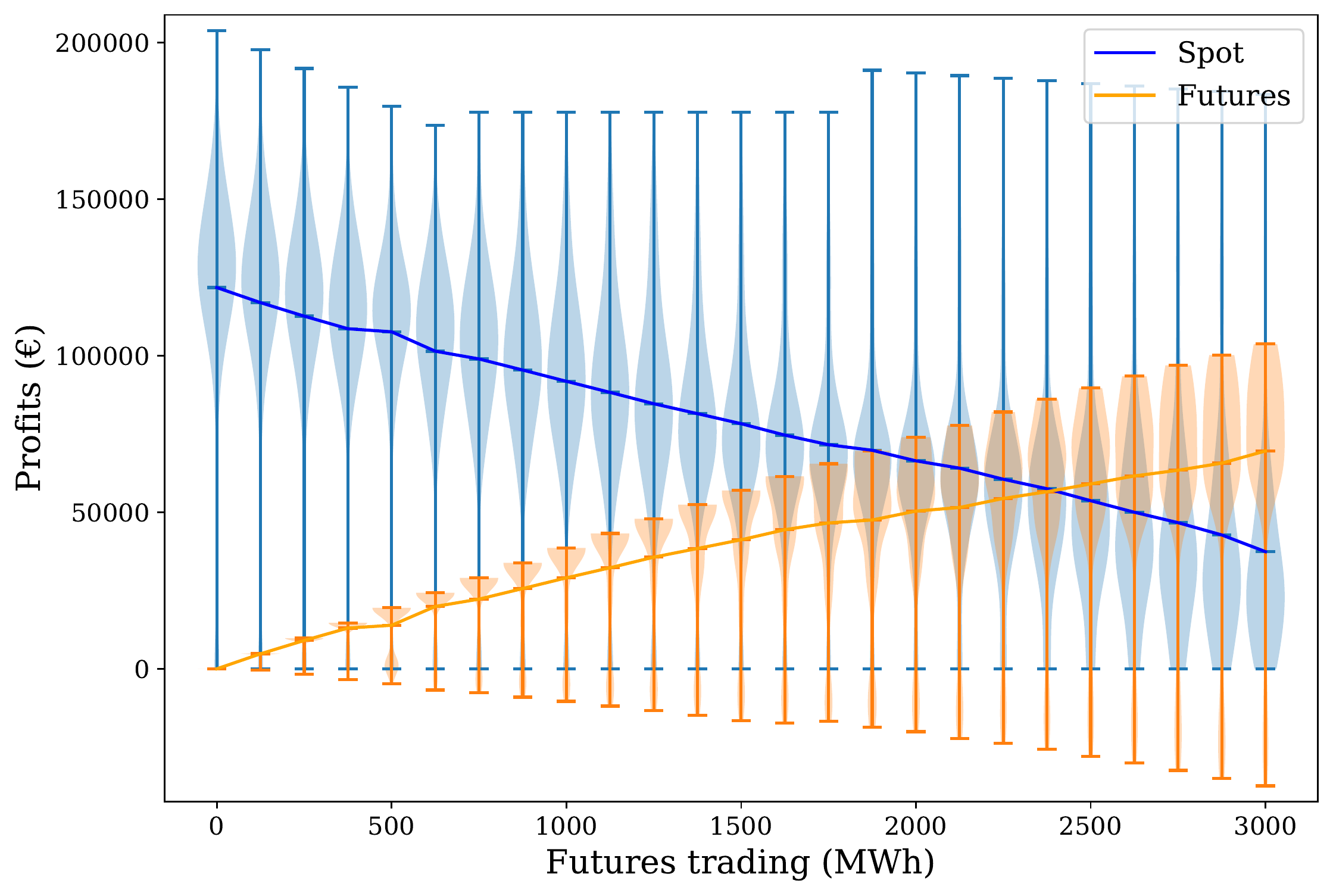}
    	\caption{Spot vs futures profit}
    	\label{fig:profit_dist_spot_futures}
	\end{subfigure}%
	\caption{Profit distributions}\label{Fig_:profit_tot}
\end{figure}

\begin{figure}[H]
	\centering
		\includegraphics[height=7cm]{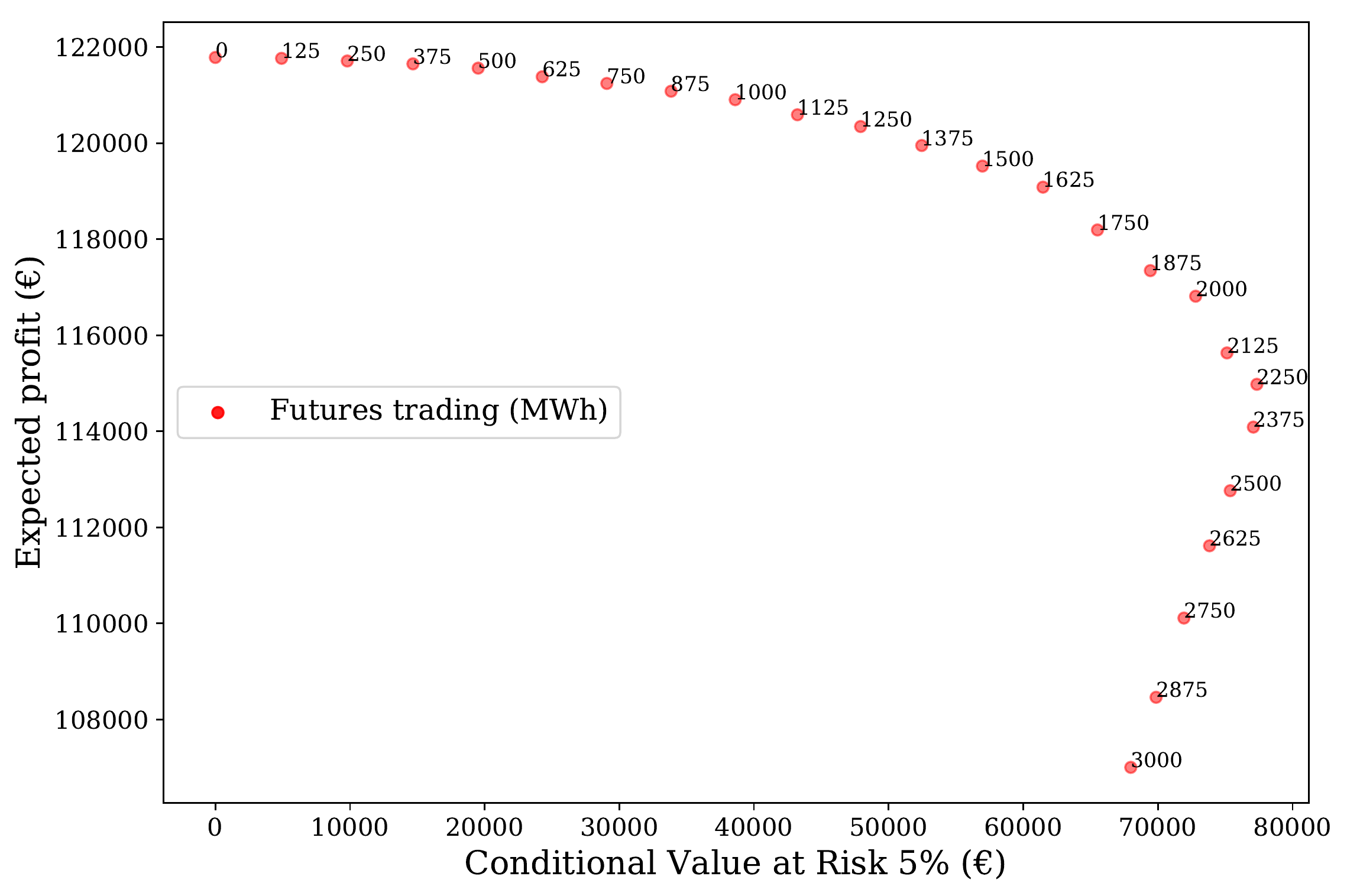}
\caption{Efficient Frontier}
	\label{fig:effic_front}
\end{figure}

This result can be
immediately observed in Figs. \ref{Fig_:profit_tot} and \ref{fig:effic_front}. In particular, Fig. \ref{fig:profit_dist} presents the profit distribution (through violin plots) for each value of $q^{ft}$ together with their mean value (blue line) and its 5\% CVaR (red line). Fig. \ref{fig:profit_dist_spot_futures} split this profit into the spot market (blue) and the futures market (orange) ones. While Fig. \ref{fig:effic_front} depicts the corresponding profit expected vs the CVaR values for each level of futures trading.  

Increasing the futures commitment
from zero to $2250$, we see that the CVaR of profits regularly increase as
well (Figs. \ref{fig:profit_dist} and \ref{fig:effic_front}). Such an outcome is not in general true in this kind of problems, as it
has been observed in the literature. Indeed, as we can observe in the same
figures, also our results show that there is a rational maximum level for
the number of futures, and exceeding it generates an increase of risk (Fig. \ref{fig:effic_front}). We
point out that to correctly interpret those two figures, as well as all
others showing the CVaR measure in this paper, this variable looks at the
left tail of the profit distribution. As such, rational risk preference is
for increasing CVaR values (actually our CVaR measure should be better
referred as a safety measure, rather than a risk one). Indeed Fig. \ref{fig:profit_dist_spot_futures} shows how the increase of the profit obtained in the futures market ``improves'' the left tail of the total profit distribution (Fig. \ref{fig:profit_dist}).

\begin{figure}[H]
	\centering
	\includegraphics[height=7cm]{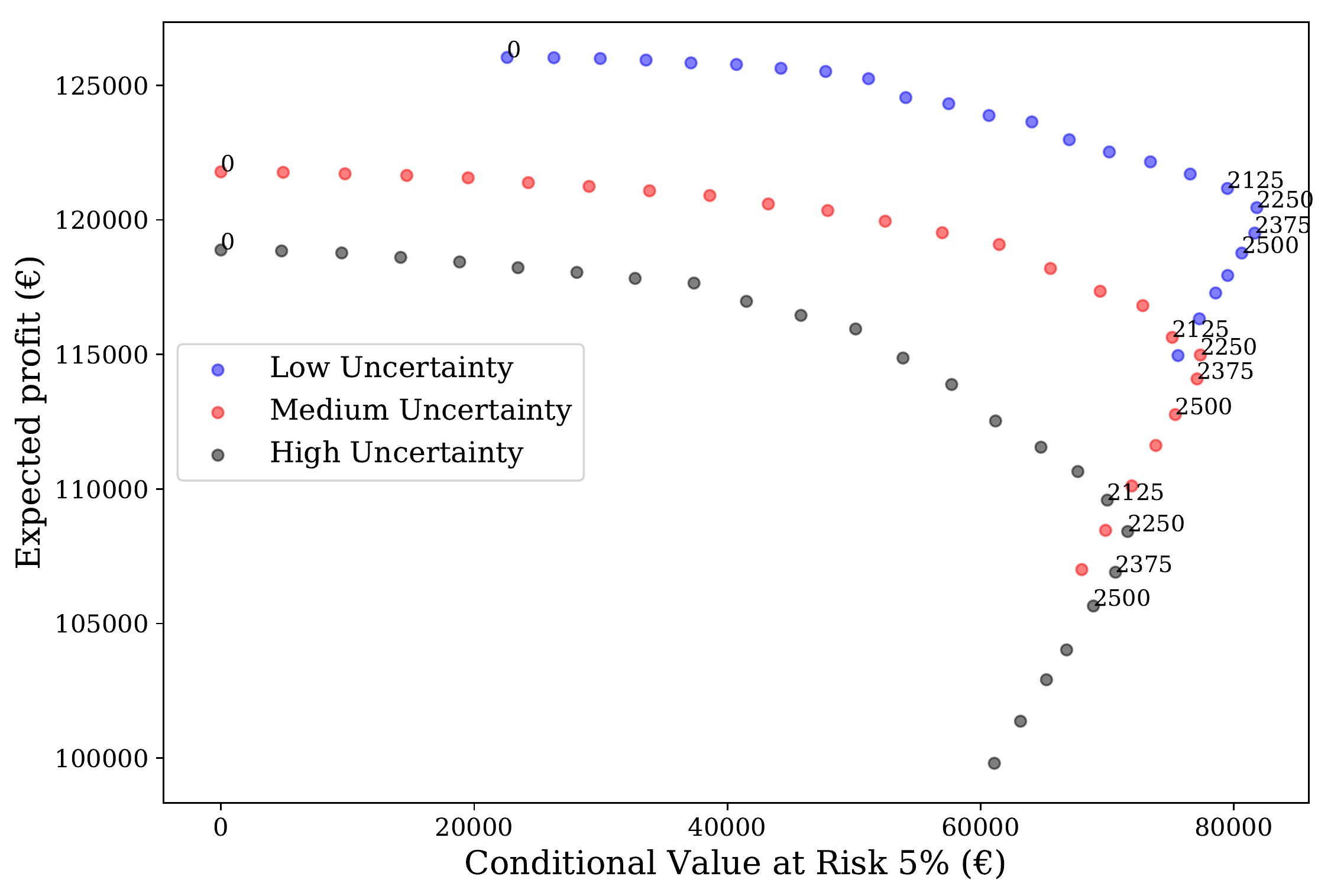}
	\caption{Efficient frontier for different demand uncertainty levels}
	\label{fig:EF_demand}
\end{figure}

The main reason why futures reduce risk seems to depend on the good
protection from the volatility of the demand. In particular, Fig. \ref{fig:EF_demand} compares the relationship expected vs profit CVaR when we modify the level of demand uncertainty (demand scenarios for these cases are generated with and increase/decreased of 20\% of their variability).

Indeed, uncertain demand worsens both profits expectation and risk. 
The worsening of risk is not surprising. However, why increasing volatility of demand reduces
expected profits is a more involved question. Its understanding explains why
futures are effective in hedging this source of risk. A starting point can
be that the (optimized) marginal technology on the spot market decided by
the strategic producer tends to be more expensive in general than the
average cost technology. In this case, the impact of higher demand
uncertainty can be asymmetric: price increase (driven by high demand levels)
would be out-weighted in size by price drops (driven by low levels of
demand), which easily fall to zero. The overall impact of higher demand
uncertainty is therefore to depress expected profits on the spot market. On
the contrary, the futures component of profit remains unaffected by the
volatility of the demand: futures take delivery priority in the case of low
demand and enjoy fixed selling price. So, they brilliantly avoid the zero
profit scenarios (see the red line of CVaR raising almost linearly from $0$
in Fig. \ref{fig:profit_dist}).

\begin{figure}[H]
	\centering
	\includegraphics[height=7cm]{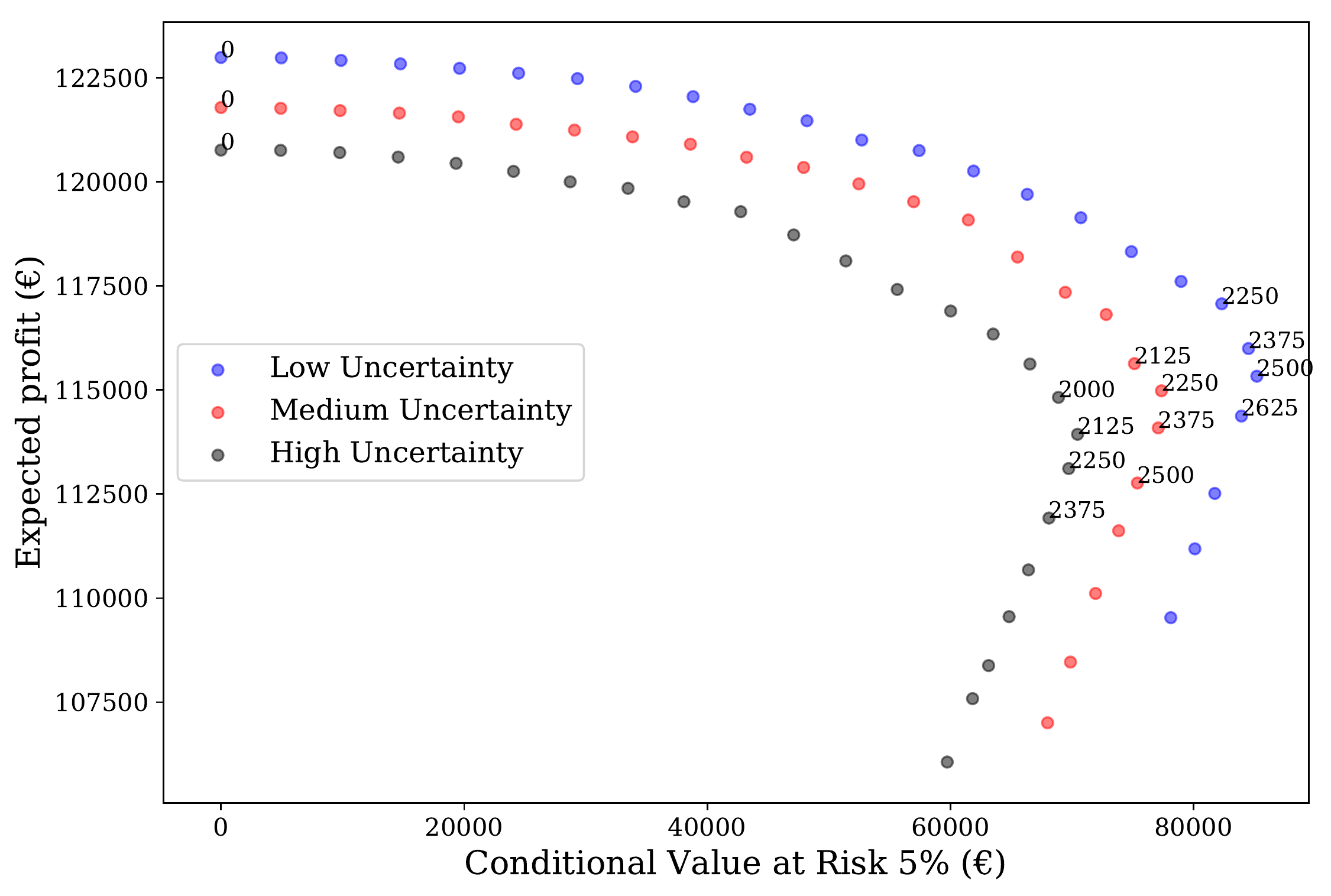}
	\caption{Efficient frontier for different renewable generation uncertainty levels}
	\label{fig:EF_renew}
\end{figure}

The volatility of renewables output worsens both expected profit and risk
and futures provide an intermediate protection against it. This is observed in Fig. \ref{fig:EF_renew} which compares the trade-off expected vs profit CVaR for different levels of renewables uncertainty (scenarios for renewable capacity are generated with and increase/decreased of 20\% of their variability). 

If the futures commitment cannot be satisfied with low-cost renewable plants
(which is a general rule to keep the spot price as high as possible), the
expensive ones become necessary, so the profits resulting from futures are
directly affected by the volatility of the renewables. However, the spot
market channel suffers this source of risk even more. Indeed, spot profits
tend to collapse both when renewables are very low and when they are very
high, while futures suffer only in the case of low renewables output: in the
first case because the producer loses the highest unit profit source; in the
second case, because large renewables outputs depress price. Overall, if the
output of renewables gets very volatile both spot and futures channels are
negatively affected. However, fixed price contracts hedge against the
positive shocks of such variable.

\begin{figure}[H]
	\centering
	\includegraphics[height=7cm]{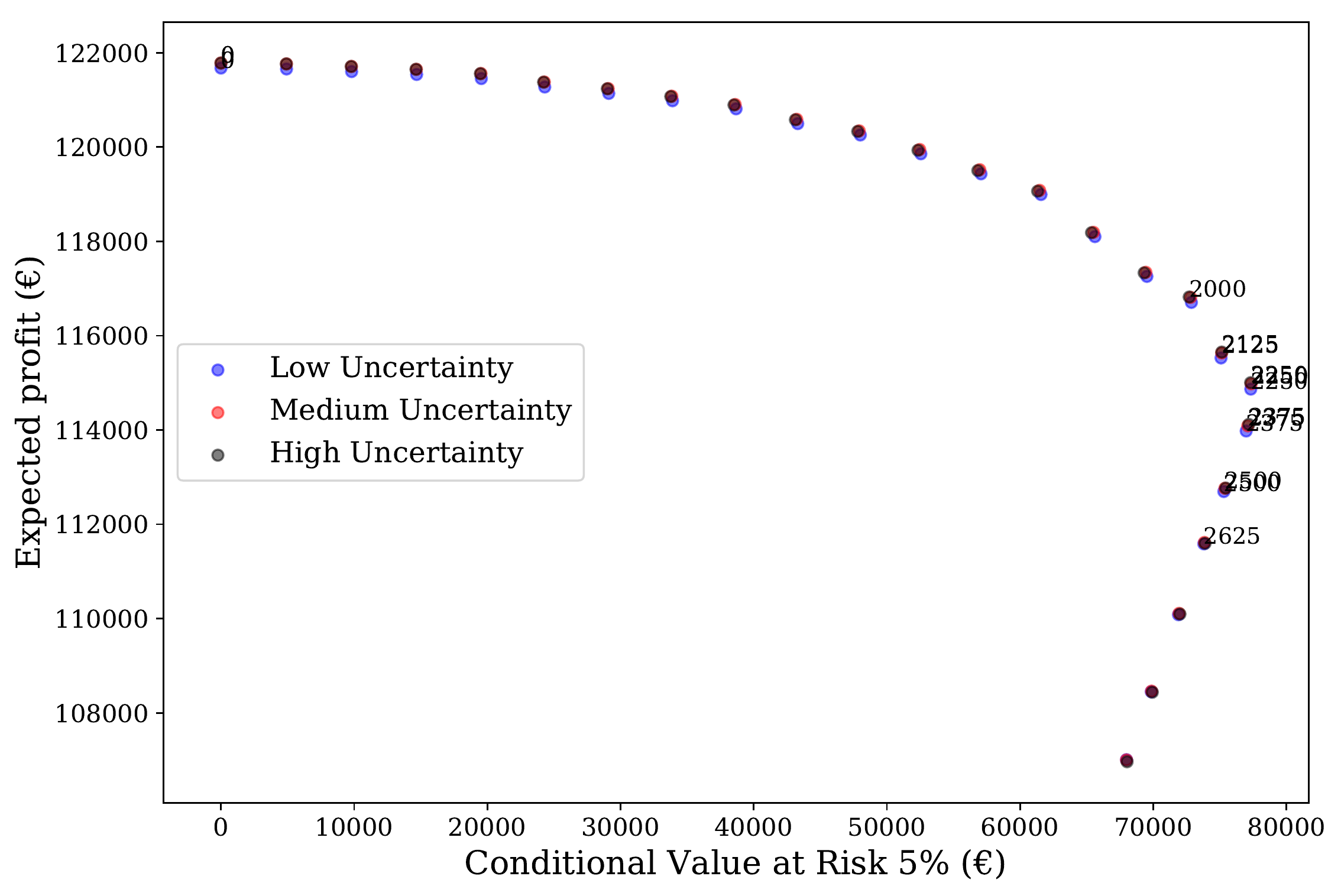}
	\caption{Efficient frontier for different cost uncertainty levels}
	\label{fig:EF_cost}
\end{figure}

Finally, futures do not provide any protection from the volatility of direct
generation costs. Vice versa, spot market channel provides an efficient
protection against this risk, since electricity demand is strictly inelastic
in the short run and producers can easily transfer to the selling price all
the shocks on the direct costs. To neutralize generation costs volatility
also on the futures channel, our producer gives priority to low costs
technologies (renewables and nuclear) as they are not subject to the
volatility of direct generation costs. The impressive results of such
optimal combination of sales mix and generation plan can be observed in Fig.
\ref{fig:EF_cost}, where the different levels of volatility of the direct costs do not
modify the optimal pattern of the benchmark case (Medium Uncertainty).

\subsubsection*{\textbf{Futures reduce expected profits}}
This result was only partly
expected. Indeed, and some circumstances, which do not apply to the model
studied here, futures can be used with a speculative intent, to raise the
expected profit of a power producer. This is in particular the case when the
futures price of electricity show positive risk premium. However, as we discuss in the
following, the hypothesis adopted in this model, namely the ``naive
expectation" of competitors, prevent such market situation. The explanation
of the negative impact of futures on the expected profits of the producer is
therefore analytically clear. Since the futures price is regularly lower
than the expected spot price, every unit of energy sold under the former
channel produces a lower profit (in expectation) than if it was sold through
the second channel (see the columns of the futures and the expected spot
prices in Table \ref{summarySM}). Overall, the futures contracts under the settings of
this model provide a classical profit / risk trade-off.

\begin{table}[ht]\caption{Solution Summary}\label{summarySM}
	\centering
	\resizebox{\columnwidth}{!}{%
	\begin{tabular}{ |+c |^c | ^c | ^c | ^c | ^c | ^c | ^c | ^c | ^c | ^c |} 
		\hline
		Futures & Expected & CVaR &  Futures & Exp. Spot  & \multicolumn{2}{c|}{Nuclear Prod.} & \multicolumn{2}{c|}{Renew. Prod.} & \multicolumn{2}{c|}{Conven. Prod.} \\
		\cline{6-11}
		 trading & Profit & Profit & Price & Price & Spot & Fut.  & Spot & Fut. & Spot & Fut.  \\ 
		 (MWh) & (\euro) & (\euro) & (\euro/MWh) & (\euro/MWh) & (MWh) & (MWh)  & (MWh) & (MWh) & (MWh)& (MWh)  \\ \hline \hline 
\rowstyle{\bfseries}0 & 121782.99 & 0.00 & 39.49 & 39.49 & 1238.40 & 0.00 & 1735.33 & 0.00 & 1219.96 & 0.00  \\ 
250 & 121709.62 & 9798.56 & 39.20 & 39.49 & 1022.99 & 208.92 & 1692.51 & 37.00 & 1215.88 & 4.08  \\ 
500 & 121560.37 & 19520.86 & 39.04 & 39.49 & 1106.25 & 129.16 & 1442.51 & 284.81 & 1133.93 & 86.03  \\ 
750 & 121241.92 & 29074.50 & 38.77 & 39.49 & 1007.91 & 226.63 & 1504.34 & 226.30 & 922.89 & 297.08  \\ 
1000 & 120903.82 & 38608.48 & 38.61 & 39.49 & 1041.27 & 197.18 & 1293.13 & 440.19 & 857.34 & 362.63  \\ 
1250 & 120344.84 & 47921.58 & 38.34 & 39.49 & 945.87 & 295.85 & 1249.09 & 478.74 & 744.55 & 475.41  \\ 
1500 & 119521.32 & 56970.15 & 37.98 & 39.49 & 1018.75 & 217.97 & 966.41 & 762.45 & 700.38 & 519.58  \\ 
1750 & 118191.92 & 65513.23 & 37.44 & 39.49 & 982.74 & 257.31 & 789.92 & 931.66 & 658.94 & 561.02  \\ 
2000 & 116811.79 & 72804.80 & 37.00 & 39.49 & 934.28 & 303.50 & 622.95 & 1105.13 & 628.59 & 591.37  \\ 
\rowstyle{\bfseries}2250 & 114978.26 & 77350.62 & 36.46 & 39.49 & 920.90 & 313.36 & 463.69 & 1277.24 & 560.56 & 659.40  \\ 
2500 & 112763.29 & 75374.22 & 35.88 & 39.49 & 868.57 & 363.00 & 274.65 & 1469.62 & 553.16 & 667.38  \\ 
2750 & 110112.18 & 71935.13 & 35.25 & 39.49 & 749.86 & 491.31 & 169.84 & 1573.73 & 537.49 & 684.97  \\ 
3000 & 107003.69 & 67986.60 & 34.59 & 39.37 & 575.54 & 667.53 & 164.25 & 1586.43 & 484.25 & 746.04  \\ 
		\hline
	\end{tabular}
}
\end{table}

\subsubsection*{\textbf{Expected spot market prices are constant with respect to the
commitment to futures}}
More precisely, the producer, whatever the sales
mix choice, succeeds to fix the same price on the spot market (see column 5
of expected spot prices in Table \ref{summarySM}). To obtain such result it suffices to
him that, whatever is the amount of the futures commitment, the spot price
is fixed by the same technology, namely the right-most possible one on the
merit order given the level of the demand. This is indeed possible to him
since any commitment to futures corresponds exactly to an equal reduction
of the residual demand on the spot market. This leaves the equilibrium point
on the spot market unchanged.

\subsubsection*{\textbf{Futures prices / market show negative risk premium}}
\begin{figure}[H]
	\begin{subfigure}[b]{0.5\textwidth}
		\includegraphics[height=5.5cm]{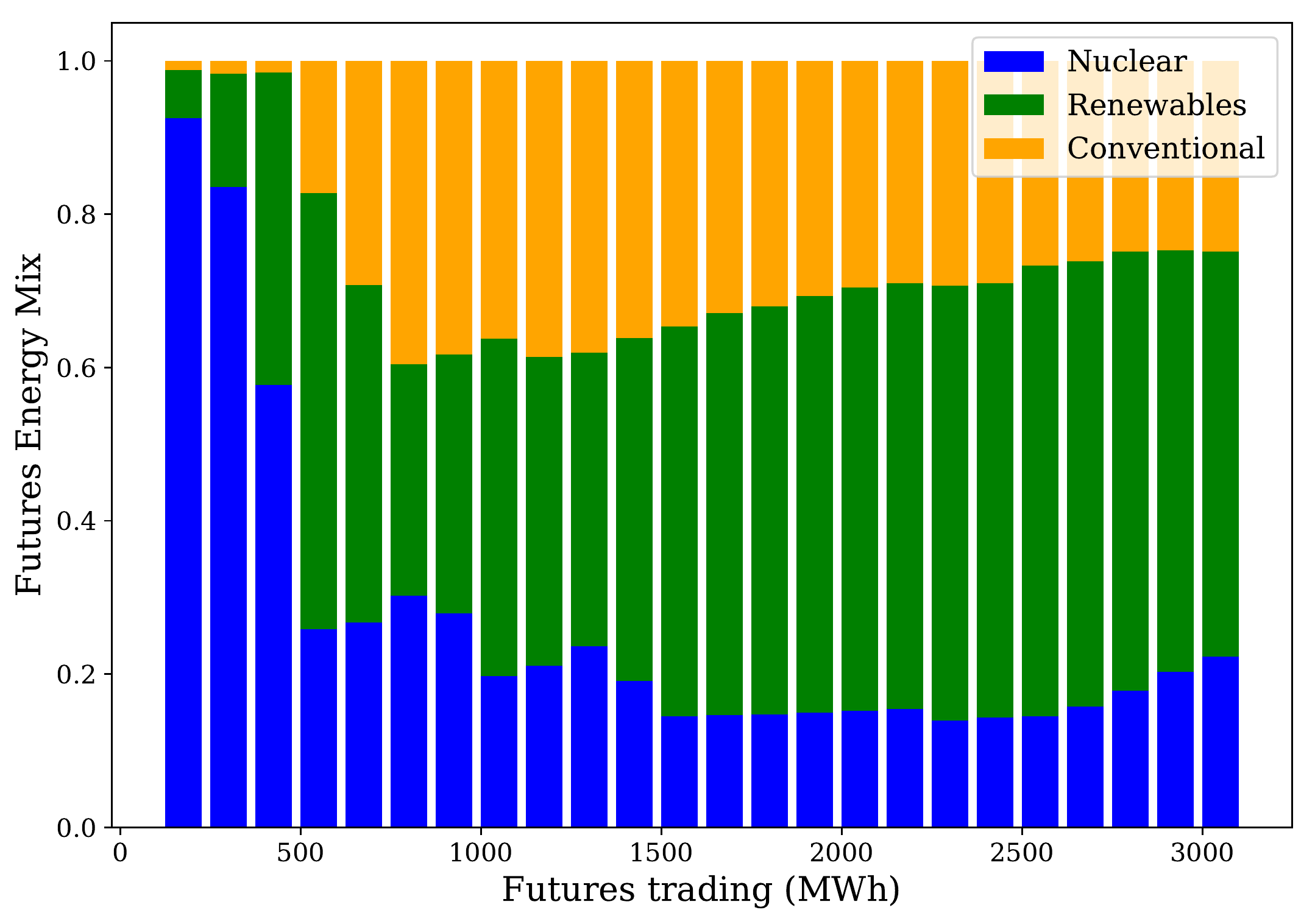}
		\caption{Futures market}
		\label{fig:fut_mix}
	\end{subfigure}
	\hskip 0ex
	\begin{subfigure}[b]{0.5\textwidth}
		\includegraphics[height=5.5cm]{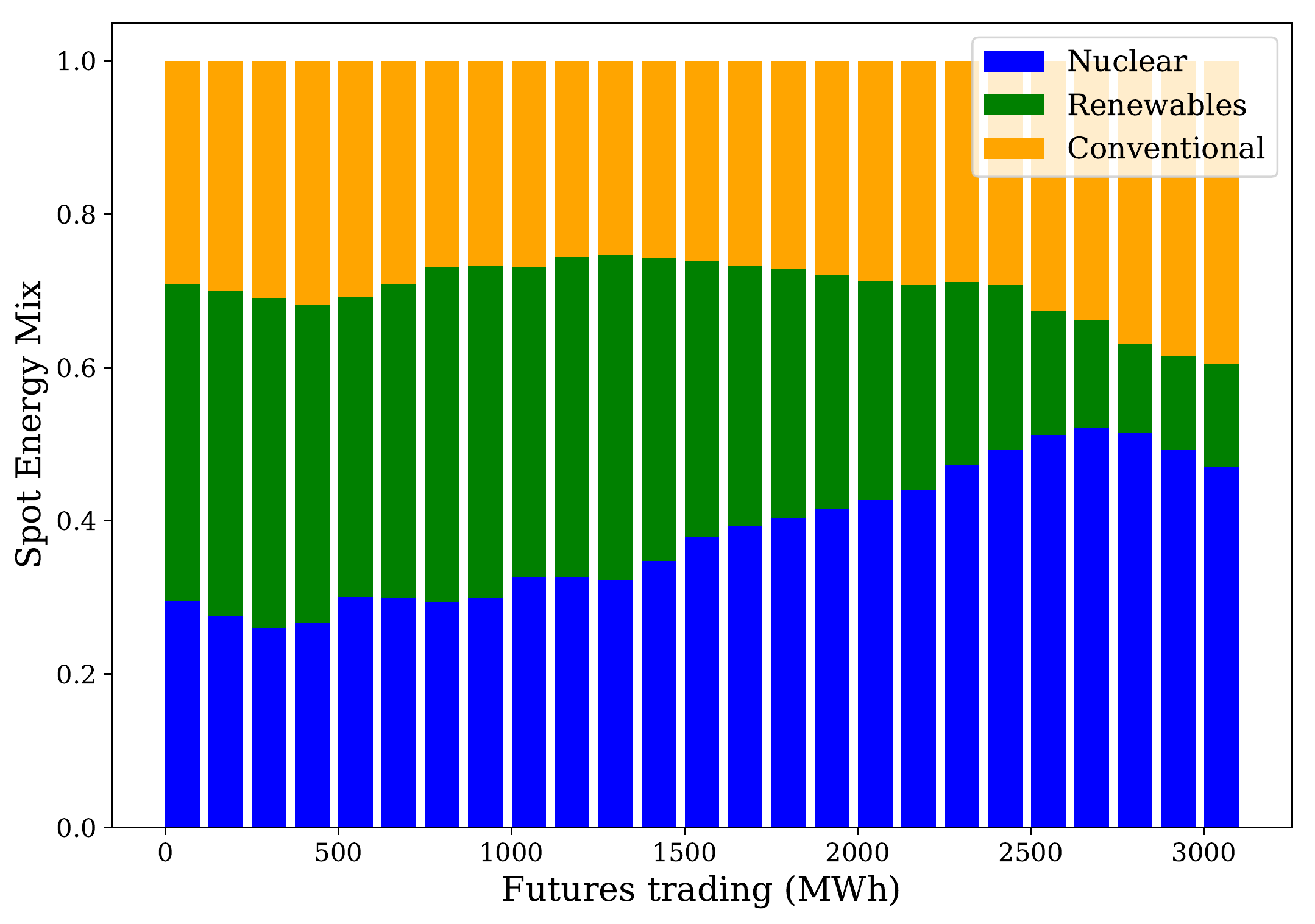}
		\caption{Spot market}
		\label{fig:spot_mix}
	\end{subfigure}%
	\caption{Technology production mix: Strategic Model}\label{fig:energy_mix}
\end{figure}

\begin{figure}[H]
	\begin{subfigure}[b]{0.5\textwidth}
		\includegraphics[height=5.5cm]{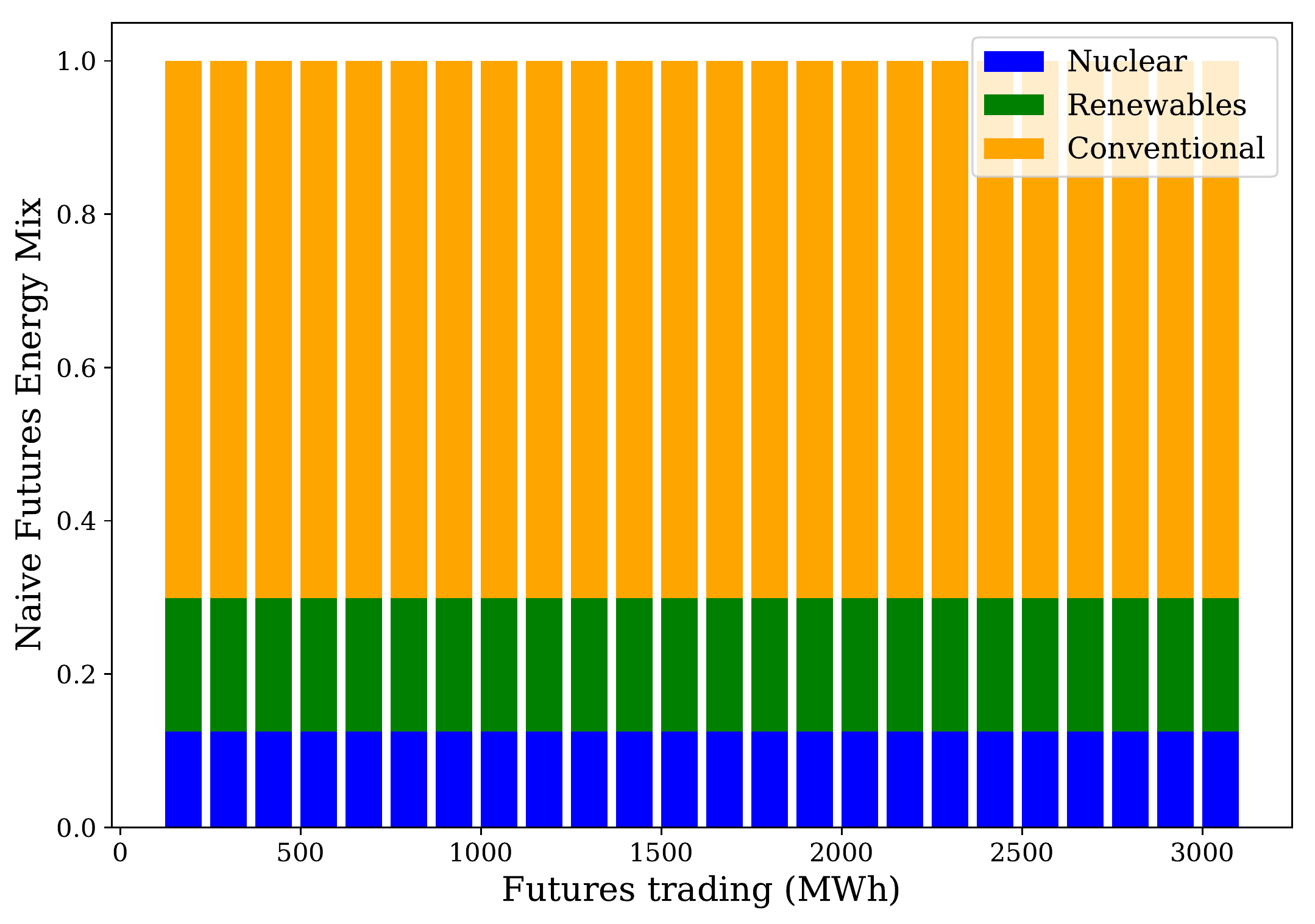}
		\caption{Futures market}
		\label{fig:fut_mix_naive}
	\end{subfigure}
	\hskip 0ex
	\begin{subfigure}[b]{0.5\textwidth}
		\includegraphics[height=5.5cm]{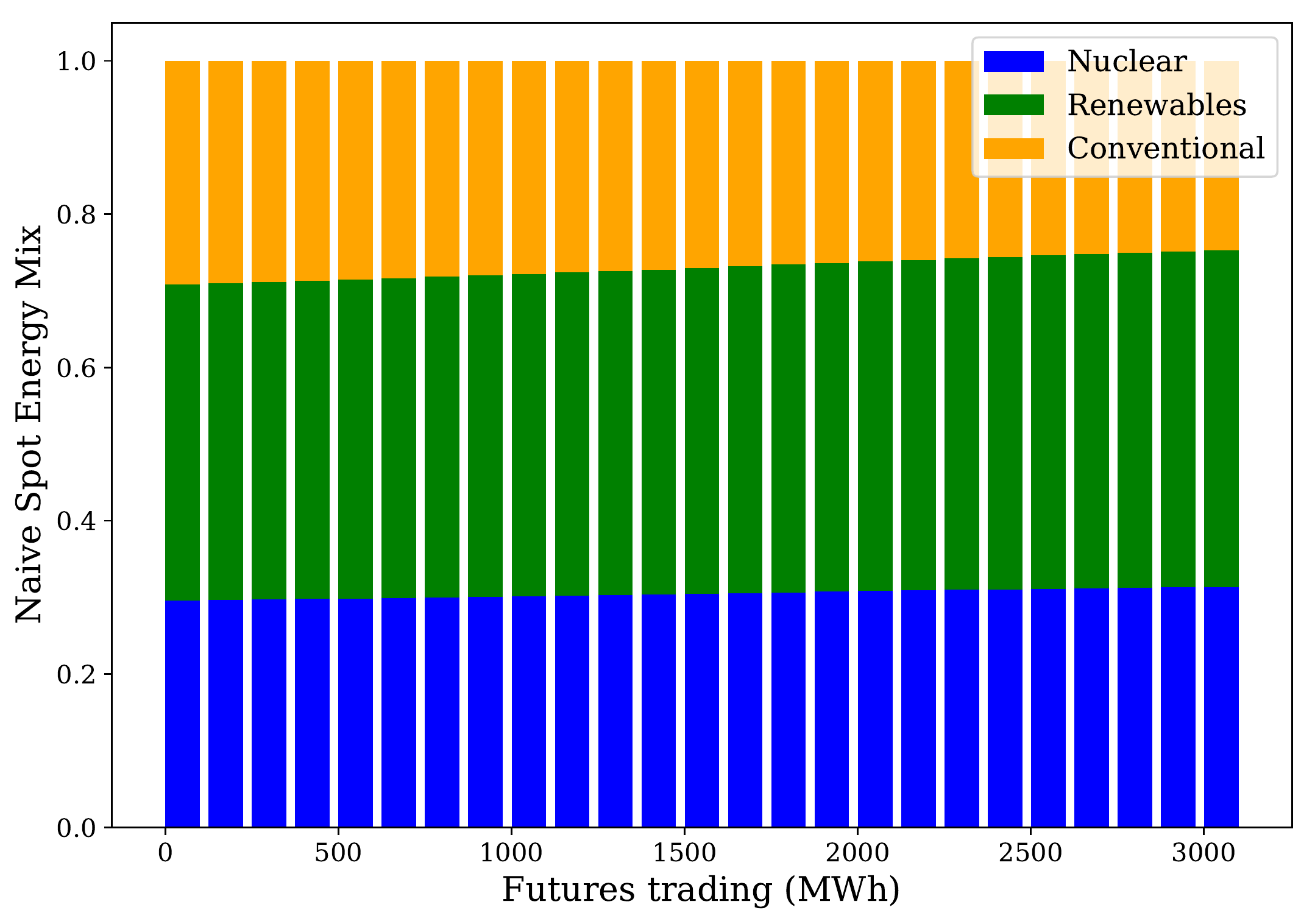}
		\caption{Spot market}
		\label{fig:spot_mix_naive}
	\end{subfigure}%
	\caption{Technology production mix: Naive Model}\label{fig:energy_mix_naive}
\end{figure}

This result (see again the
columns of the futures and the expected spot prices in Table \ref{summarySM})
can be explained in the first place by the so called ``naive expectation"
hypothesis that is assumed here for the competitors. Under such hypothesis it is assumed that competitors do not correctly anticipate the  generation plan of the strategic producer. In particular, given their information gap, they end up to underestimate the use of low-cost technologies by the strategic producer to satisfy his commitment to futures. This can be observed comparing Figs. \ref{fig:energy_mix} and \ref{fig:energy_mix_naive}, which represent respectively the optimal generation plan adopted by the strategic producer and that assumed by the competitors under the naive hypothesis.

The immediate consequence is that the competitors obtain an over estimation of low cost technologies on the spot market, which of course reduces the value of the expected spot price and ultimately that of the futures price, given by the fundamental equation (\ref{bi_cons_3}).

A partial mitigation to the doom of negative risk premium can be that played by the parameter $\alpha$ appearing in
the balance constraint (\ref{naive_bal}). This parameter interacts with the naive hypothesis letting the competitors assume higher or lower impact on the residual demand
(for the spot market) than the exact quantity of futures decided by the
strategic producer. Such bias can occur on the real markets, because of the
reaction of the competitors, which can emphasize or de-emphasize the
decision of the strategic producer. In the first case, competitors react
selling more futures than originally planned and so they reduce the residual
demand (on the spot market) even more than the strategic producer
determines. The contrary occurs in the opposite case. If such parameter is
equal to $1$, which is the value adopted throughout in this paper, the
demand is reduced exactly by the same amount of the futures commitment, if
it is higher or lower than $1$, we model the emphasizing or the
de-emphasizing reaction respectively. In particular, values of $\alpha<1$ are
expected to sustain the futures price and to re-balance negative risk premium.

While our general findings point out that futures prices show negative risk premium, the well-known paper of \cite{BL2002} supplies a model where futures prices show both positive and negative risk premia. The origin of such difference can be traced to the specific hypotheses about the producers. More precisely, the two cited authors assume that all producers are the same with respect to capacity (actually with no upper limit) and same cost function.  On the contrary, we consider producers with different generating technologies and sizes and let the strategic producer be large enough to affect the market price. Besides, the two authors take the demand as the only source of uncertainty, while we consider a larger set of uncertainty, including both generation levels and direct costs. 
Assuming that all producers are the same is not compatible with the ``naive'' hypothesis adopted in this paper. Eliminating any information gap between rivals and the strategic producer would let the latter forecast correctly the expected market price. Under such setting our model would perfectly equate futures price to expected spot price, eliminating chances for both positive and negative risk premium.
So \cite{BL2002} model is not directly comparable to our. That model grounds on the distribution properties of the demand to explain why futures price can positive or negative risk premia. Our model focuses on the information gap between the strategic producer and his/her rivals. As long as such information gap underestimates the market power of the strategic producer, our model predicts that futures prices are in backwardation. 

\subsection{Data analysis of the low-level solutions}

To add a clearer view of the output of the optimal strategy commented in the
previous points, we also add a particular analysis based on a combination of
cluster and discriminant analysis. The object of this data analysis has
been the optimal solutions adopted in each of the 300 random scenarios (at
the low-level problem), for two commitment of futures, namely 250 and 2250.
In particular, the observed variables are those listed in the first column of
Table \ref{summayCLa}.

\textbf{Cluster analysis}. The 300 scenarios have been first grouped in 4
clusters, by means of hierarchical criterion based on Ward method. Only two
variables have been used to this purpose: the random level of the demand and
of the total renewables output. These two variables, as we have seen in the
previous section, are those which have a significant impact on both the
expected profit and the CVaR of profit for our producer. Two additional
clusters have been added manually, namely clusters 5 and 6, separating them
manually from clusters 2 and 3 respectively. The reason of this intervention
was to highlight separately the particular feature of such sub-scenarios,
namely the outcome of a spot market closing at zero price.

\begin{figure}[H]
	\centering
	\includegraphics[height=7cm]{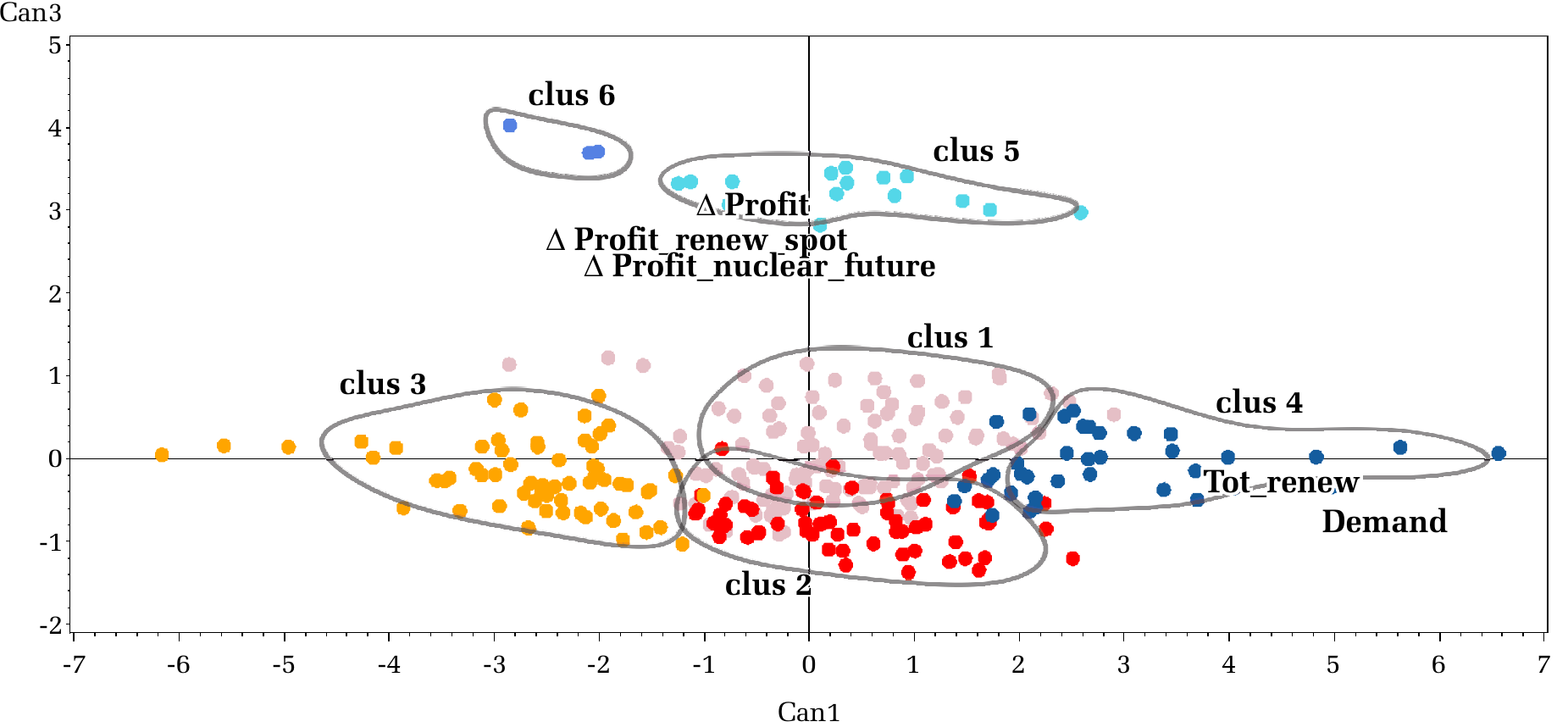}
	\caption{Discriminant analysis}
	\label{fig:clus6a}
\end{figure}

\textbf{Discriminant analysis.}, The main purpose of this analysis was to
work out unified picture of the 6 clusters on two discriminant axis. Figure (%
\ref{fig:clus6a}) provides such a representation. The horizontal axis is
positively correlated with high levels of the demand and renewables output (%
\textit{Tot\_renew} and \textit{Demand}, in the figure). The vertical one
with the impact of passing from a futures commitment of 250 to 2250 for:

\begin{itemize}
\item overall profit ($\Delta $\textit{\ Profit}),

\item spot profit for renewables ($\Delta $ \textit{Profit\_renew\_spot})

\item futures profit for nuclear plants ($\Delta $ \textit{%
Profit\_nuclear\_futures}).
\end{itemize}

This analysis confirms and highlights a few interesting facts.

\textbf{Futures contracts really make a difference in the worst-case
scenarios}. We can see this observing in particular the scenarios in cluster
5 and 6. In these scenarios the spot market closes at zero price, provides no
profit and the only income is based on the futures commitment. So, with
a short position in futures of 250 the income is 9'799'000 \euro\ and it
grows to 84'400'000 \euro\ passing to a position of 2250 (see Table \ref%
{summayCLa}), which a substantial improvement. The neatly separated position
of these clusters from all the other ones in Figure (\ref{fig:clus6a}) shows
clearly such an impact. 

\textbf{Futures contracts also help when renewables output is large}. As we
have already commented, large renewable output can produce severe drops in
price, even when demand is average. This is of course a serious problem for
an electricity producer. We can see this situation in the case of the
scenarios in cluster 2, where the renewables output averages to 8818 $MWh$
(see Table \ref{summayCLa}).  Joint with an average demand of 16'879 $MWh$
this produces a spot price of $38.9$ \euro , which is a drop of $-11.2\%$
with respect to the case of cluster $4.$ However, the drop in the overall
profit is equal to $-3^{\prime }691$ \euro\ and $-13513$ \euro , that is
only $-2.7\%$ and $-9\%\,$with respect to the levels of the profit in
cluster 4, adopting a futures position of 2250 and 250 respectively.

\textbf{The general rule is to satisfy futures commitment with low costs
technologies}. As we can see in the cases of the scenarios in clusters 2, 5
and 6 (the lowest prices scenarios) passing from 250 to 2250 short futures
position helps substantially to improve the profit performance. Observing
Table \ref{summayCLa} the optimal generation in these cases plan consists
of moving low costs technologies from spot market to futures in order to
guarantee that spot price keeps as high as possible. In the case such
strategy does not succeed to keep the spot price high enough or even zero
(i.e. as in clusters 5 and 6), the income of the futures channel provides a
good insurance anyway.

\begin{table}[tbp]\caption{Clusters average values}\label{summayCLa}
\centering%
\begin{tabular}{lrrrrrr}
\hline
& \multicolumn{1}{l}{clust1} & \multicolumn{1}{l}{clust2} & 
\multicolumn{1}{l}{clust3} & \multicolumn{1}{l}{clust4} & \multicolumn{1}{l}{
clust5} & \multicolumn{1}{l}{clust6} \\ \hline
Tot\_renew\_out ($MWh$) & 6'413 & 8'818 & 5'350 & 8'324 & 8'842 & 5'993 \\ 
Str\_Renew\_out ($MWh$) & 1'605 & 2'198 & 1'337 & 2'061 & 2'196 & 1'582 \\ 
Demand ($MWh$) & 20'770 & 16'879 & 15'135 & 24'175 & 12'435 & 10'265 \\ 
\hline
Profit - Fut.250 (\euro 000) & 129'676 & 135'941 & 110'017 & 149'454 & 9'799
& 9'799 \\ 
Profit - Fut.2250 (\euro 000) & 114'449 & 130'378 & 99'362 & 134'069 & 82'044
& 82'044 \\ 
$\Delta $ Profit (\euro 000) & -15'226 & -5'562 & -10'655 & -15'385 & 72'245
& 72'245 \\ \hline
Spot Price (\euro$/MWh$) & 43.7 & 38.9 & 41.5 & 43.8 & 0.0 & 0.0 \\ \hline
$\Delta $ Str\_Spot\_Nucl ($MWh$) & -18 & -134 & -139 & -49 & -458 & -928 \\ 
$\Delta $ Str\_Spot\_Renew ($MWh$) & -1'138 & -1'510 & -1'136 & -1'176 & 
-1'302 & -902 \\ 
$\Delta $ Str\_Spot\_Conven ($MWh$) & -845 & -356 & -724 & -775 & 0 & 0 \\ 
$\Delta $ Str\_Fut\_Nucl ($MWh$) & 18 & 134 & 139 & 49 & 476 & 735 \\ 
$\Delta $ Str\_Fut\_Renew ($MWh$) & 1'138 & 1'510 & 1'136 & 1'176 & 1'524 & 
1'265 \\ 
$\Delta $ Str\_Fut\_Conv ($MWh$) & 845 & 356 & 724 & 775 & 0 & 0 \\ \hline
$\Delta $ Profit\_nuclear\_spot (\euro 000) & -656 & -4'913 & -5'523 & -2'067
& 0 & 0 \\ 
$\Delta $ Profit\_renew\_spot (\euro 000) & -49'574 & -58'309 & -46'489 & 
-51'356 & 0 & 0 \\ 
$\Delta $ Profit\_conv\_spot (\euro 000) & -4'382 & -1'471 & -3'032 & -4'042
& 0 & 0 \\ 
$\Delta $ Profit\_nuclear\_futures (\euro 000) & 43 & 4'291 & 4'505 & 1'148
& 17'220 & 26'586 \\ 
$\Delta $ Profit\_renew\_futures (\euro 000) & 41'402 & 55'011 & 41'350 & 
42'822 & 55'025 & 45'659 \\ 
$\Delta $ Profit\_conv\_futures (\euro 000) & -2'060 & -172 & -1'466 & -1'891
& 0 & 0 \\ \hline
\end{tabular}%
\end{table}

\color{black}

\color{black}

\section{Conclusions}
    The model proposed here allows us to analyse several interesting details of the joint sales mix and generation plan management decisions, as well as to understand some significant facts on power markets.
    In this model we let futures prices be an endogenous (risk neutral) expectation of the equilibrium spot price of electricity. The way such expectation is estimated is therefore of crucial importance. As long as market players cannot fully predict the optimal management decisions of a major producer, they introduce a negative bias in spot price expectation and set a negative risk premium. To the best of our knowledge this is the first attempt to introduce an information gap hypothesis to model negative risk premia of futures in the power sector. Besides the theoretical interest on its own, this has significant consequences on the role that futures and in general fixed prices contracts can play for to a large energy producer, that is one that can at least partially control the spot price of electricity. In particular, as we observe in the paper, fixed price contracts end up providing just a risk hedging role. This is not an obvious conclusion in the case of power markets, where opportunities to partially control the price of electricity exist, since storability is limited, demand is highly rigid (especially in the short run) and in most cases these markets are highly concentrated. As the empirical and theory literature have shown, futures prices on electricity markets can show as well be positive risk premium, and in this case a strategic producer can enter short futures positions following (rationally) a speculative purpose.
    In any case, the ``naive hypothesis'' provides a new potential explanatory factor to the empirical models of futures risk premia, that can be tested empirically comparing them in different markets, with different concentration levels and during periods of different volatility and market predictability. 
    
    Our results show that the risk reduction provided increasing the short position of futures in the sales mix is substantial, with the CVaR of lowest profits improving about linearly with the level of the futures. At the same time, we show that the futures position cannot be increased indefinitely, since there is a maximum threshold which is not optimal to exceed. Such threshold depends on several fundamental market and idiosyncratic variables (i.e., of the strategic producer), as well as on the volatility of the renewables output.
    We show the mechanics through which futures deploy an effective hedge of risk. As it can be expected, futures reduce the expected profit in the scenarios where profit level is high but guarantee substantial profit levels in the lean scenarios. The interesting point highlighted here is how the futures positions combine optimally with the generation plan. In short, the strategic producer satisfies the commitment in futures giving priority to low generation costs technologies. In this way he guarantees a safe profit in the lean scenarios and, at the same time, maximizes the expectation of the spot price of electricity.
    The analysis also pointed out the different contributions of the sales mix to reduce the CVaR of profits. In particular, we considered three (joint) major sources of uncertainty in this paper: volatility of demand, of renewables output and direct generation costs (i.e., the cost of fuels). We showed that taking short futures position strongly improves the CVaR in the scenarios of low demand. Besides, futures offer a partial risk reduction with respect to the volatility of renewables output, since they effectively protect only the cases of high output levels. On the contrary, futures offer no protection against the volatility of generation costs.

\end{document}